\documentclass[11pt]{amsart}
\usepackage{hyperref}
\usepackage[margin=1in]{geometry}
 \pdfoutput=1 
\usepackage{amssymb}
\usepackage{epstopdf}
\usepackage{float}
\usepackage{amsmath}

\usepackage{tikz}
\usetikzlibrary{decorations.pathreplacing,calligraphy}
\usetikzlibrary{math}
\usetikzlibrary{positioning,calc}
\usetikzlibrary{arrows}
\usetikzlibrary{decorations.markings}

\newcommand\tileA[5]{
  \draw[#5](#1,#2) rectangle (#1+#3,#2+#3);
  \begin{scope}[shift = {(#1+#3/2,#2+#3/2)}]
    \begin{scope}[rotate = {#4}]
    \draw (-#3/2,0) arc (-90:0:#3/2);
\draw (0,-#3/2) arc (180:90:#3/2);    
    \end{scope}
      \end{scope}
  }

\input{001_3.tex}
\newcommand\cupA[6]{
  \begin{scope}[shift={(#1,#2)}]
  \begin{scope}[rotate={#3},scale={#6}]    
\draw[#4,line width = #5](0,0)--(-5,0)--(-5,-5)--(-10,-5)--(-10,0);
  \end{scope}
\end{scope}      
}

\newcommand\comment[1]{}

\newcommand\tileB[5]{
  \draw[#5](#1,#2) rectangle (#1+#3,#2+#3);
  \begin{scope}[shift = {(#1+#3/2,#2+#3/2)}]
\begin{scope}[rotate = {#4}]
    \draw (-#3/2,0) -- (0,#3/2);
\draw (0,-#3/2)  -- (#3/2,0);
\end{scope}
  \end{scope}
}

\newcommand\tileC[5]{
  \draw[#5](#1,#2) rectangle (#1+#3,#2+#3);
  \begin{scope}[shift = {(#1+#3/2,#2+#3/2)}]
\begin{scope}[rotate = {#4}]
    \draw[rounded corners] (-#3/2,0) -- (0,0)-- (0,#3/2);
\draw[rounded corners] (0,-#3/2)  --(0,0) -- (#3/2,0);
\end{scope}
  \end{scope}
}

\newcommand\tileCC[5]{
  \draw[#5](#1,#2) rectangle (#1+#3,#2+#3);
  \begin{scope}[shift = {(#1+#3/2,#2+#3/2)}]
\begin{scope}[rotate = {#4}]
    \draw[rounded corners, gray, opacity=0.5, line width=0.5mm] (-#3/2,0) -- (0,0)-- (0,#3/2);
\draw[rounded corners,black, line width=1.5mm] (0,-#3/2)  --(0,0) -- (#3/2,0);
\end{scope}
  \end{scope}
}

\newcommand\tileD[5]{
  \draw[#5](#1,#2) rectangle (#1+#3,#2+#3);
  \begin{scope}[shift = {(#1+#3/2,#2+#3/2)}]
\begin{scope}[rotate = {#4}]
\draw[black, line width=1.5mm] (0,-#3)  --(0,0);
\end{scope}
  \end{scope}
  }

\def\mylist{{0}{90}{180}{270}}
\pgfmathdeclarerandomlist{myang}{\mylist}

\newcommand\twistA[1]{
  \draw(0,0) rectangle (1,1);
         \clip (0,0) rectangle (1,1);
          \begin{scope}[shift = {(0.5,0.5)}]
          \begin{scope}[scale = {0.7071},rotate={45}]
          \begin{scope}[shift = {(-0.5,-0.5)}]
            \tileA{0}{0}{1}{#1}{gray!50!white};
            \tileA{1}{0}{1}{90}{gray!50!white};
            \tileA{0}{1}{1}{0}{gray!50!white};
            \tileA{0}{-1}{1}{0}{gray!50!white};
            \tileA{-1}{0}{1}{90}{gray!50!white};            
         \end{scope}
          \end{scope}
          \end{scope}          
}

\newcommand\twistB[1]{
  \draw(0,0) rectangle (1,1);
         \clip (0,0) rectangle (1,1);
          \begin{scope}[shift = {(0.5,0.5)}]
          \begin{scope}[scale = {0.7071},rotate={45}]
          \begin{scope}[shift = {(-0.5,-0.5)}]
            \tileA{0}{0}{1}{#1}{gray!50!white};
            \tileA{1}{0}{1}{0}{gray!50!white};
            \tileA{0}{1}{1}{90}{gray!50!white};
            \tileA{0}{-1}{1}{90}{gray!50!white};
            \tileA{-1}{0}{1}{0}{gray!50!white};            
         \end{scope}
          \end{scope}
                  \end{scope}
          }

\begin{document}

\title{Hinged Truchet tiling fractals}

  \author{{
H. A. Verrill
\date{\today}
H.A.Verrill@warwick.ac.uk
    } 
    }

    \maketitle

    \section*{Abstract}
    This article describes a new
    method of producing space filling fractal curves
    based on a hinged tiling procedure.  The fractals produced can be generated
    by a simple L-system.  The construction as a hinged tiling has the
    advantage of automatically implying that the fractiles produced tessellate,
    and that the Heighway fractal dragon curve, and all the other
    curves constructed, do not cross themselves.
    This also gives a new
    limiting procedure to apply to certain Truchet tilings.
    I include the computation of the fractal dimension of one of the curves,
    and describe an algorithm for computing the sim value of the fractal
    boundary of these curves.  The curves considered have previously
    been described by \cite{Tab}, but the hinged tiling approach is new,
    as is the algorithm for computing the sim value.
    \begin{figure}
      \end{figure}
        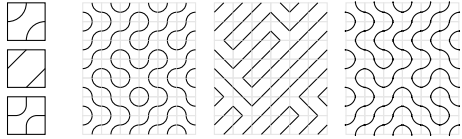
\begin{figure}[H]
      \begin{tikzpicture}[scale={0.5}]
        \tileA{0}{3.5}{1}{0}{black}
        \tileB{0}{2.25}{1}{0}{black}
        \tileC{0}{1.}{1}{0}{black}
\begin{scope}[scale = {0.5}, shift = {(4,2)}]
  \foreach \i in {0,...,5}{
    \foreach \j in {0,...,6}  {
      \pgfmathrandomitem{\myang}{myang}
      \tileA{\i}{\j}{1}{\myang}{gray!20!white}
  }}
\end{scope}

\begin{scope}[scale = {0.5}, shift = {(11,2)}]
  \foreach \i in {0,...,5}{
    \foreach \j in {0,...,6}  {
      \pgfmathrandomitem{\myang}{myang}
      \tileB{\i}{\j}{1}{\myang}{gray!20!white}
  }}
\end{scope}

\begin{scope}[scale = {0.5}, shift = {(18,2)}]
  \foreach \i in {0,...,5}{
    \foreach \j in {0,...,6}  {
      \pgfmathrandomitem{\myang}{myang}
      \tileC{\i}{\j}{1}{\myang}{gray!20!white}
  }}
  \end{scope}
        \end{tikzpicture}
      \caption{Some Truchet tiles and examples of their tilings}
      \label{troucheexamples}
\end{figure}
    
\section{Introduction}
    Space filling curves are a concept introduced by Peano in 1890 \cite{Peano},
     simplified by the Hilbert curve in 1891
    \cite{Hilbert}.
    Since then many other space filling curves have been found see e.g.,
    \cite{Sagan}.
    There are many ways to modify Hilbert's original construction e.g.,
    \cite{Ozkaraca}, and many other constructions exist.
    Space  filling curves have applications in computer science e.g.
    \cite{Bader},
    \cite{Bohm}.

    In this article, I start with a Truchet tiling,
    and iteratively produce a related tiling of the same form, in a continuous way.
    By taking the limit, a space filling curve is achieved.
    
    A {\it Truchet tiling} is 
    a set of identical square
    tiles arranged together
    in different rotations.  
    Named after Truchet who wrote about them in    1704,
    they were popularised more recently
    by \cite{Smith}, especially the top left
    tile Figure~\ref{troucheexamples}.
    These tilings are widely used in
    generative art work and graphic design, for example \cite{Krawczyk},
    \cite{Carlson},
    and also have
    been applied in computer graphics \cite{Browne}.

    In this article I restrict to the tiles
    of Figure~\ref{troucheexamples}.
    These are deformations of each other,
    and have the same symmetry groups. There is not a great
    deal of difference in the mathematics of the
    curves produced.  I will switch between tiles,
    but everything I say about one
    case applies to the others.
    
    Hinged tilings are popular and well known in
    recreational mathematics
    \cite{Frederickson},
    and have applications in origami \cite{Lang}, \cite{Barreto}.
    We use the simplest hinged tiling, that of squares, Figure~\ref{hinge}.

      \newcommand\tileAA[8]{
        \begin{scope}[shift = {(#1+#3/2,#2+#3/2)}]
\begin{scope}[scale={#7}]
          \begin{scope}[rotate={#6}]
          \draw[gray](-#3/2,-#3/2) rectangle (#3/2,#3/2);
    \begin{scope}[rotate = {#4}]
    \draw[line width = #8] (-#3/2,0) arc (-90:0:#3/2);
\draw[line width = #8] (0,-#3/2) arc (180:90:#3/2);    
    \end{scope}
    \end{scope}    
        \end{scope}
      \end{scope}
      }

            \newcommand\tileAAA[8]{
        \begin{scope}[shift = {(#1+#3/2,#2+#3/2)}]
\begin{scope}[scale={#7}]
          \begin{scope}[rotate={#6}]
          \draw[gray](-#3/2,-#3/2) rectangle (#3/2,#3/2);
    \begin{scope}[rotate = {#4}]
    \draw[line width = #8, gray!50!white] (-#3/2,0) arc (-90:0:#3/2);
\draw[line width = #8,gray!50!white] (0,-#3/2) arc (180:90:#3/2);    
    \end{scope}
    \end{scope}    
        \end{scope}
      \end{scope}
  }

    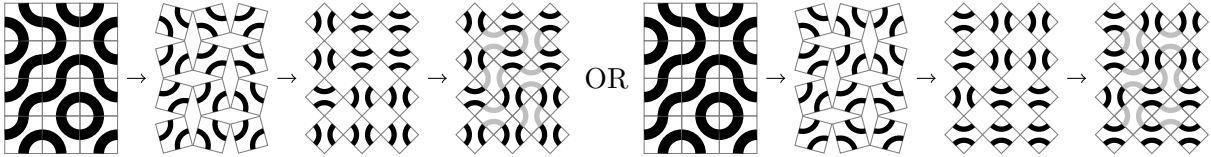
\begin{figure}
      \begin{tikzpicture}[scale = {0.5}]
        \tileAA{0}{0}{1}{0}{black}{0}{1}{4};
        \tileAA{0}{1}{1}{0}{black}{0}{1}{4};
        \tileAA{0}{2}{1}{0}{black}{0}{1}{4};
        \tileAA{0}{3}{1}{90}{black}{0}{1}{4};
        \tileAA{1}{0}{1}{90}{black}{0}{1}{4};
        \tileAA{1}{1}{1}{0}{black}{0}{1}{4};
        \tileAA{1}{2}{1}{0}{black}{0}{1}{4};
        \tileAA{1}{3}{1}{90}{black}{0}{1}{4};
        \tileAA{2}{0}{1}{0}{black}{0}{1}{4};
        \tileAA{2}{1}{1}{90}{black}{0}{1}{4};
        \tileAA{2}{2}{1}{90}{black}{0}{1}{4};
        \tileAA{2}{3}{1}{0}{black}{0}{1}{4};

\begin{scope}[shift = {(4,0)}]
        \tileAA{0}{0}{1}{0}{black}{15}{0.8}{2};
        \tileAA{0}{1}{1}{0}{black}{-15}{0.8}{2};
        \tileAA{0}{2}{1}{0}{black}{15}{0.8}{2};
        \tileAA{0}{3}{1}{90}{black}{-15}{0.8}{2};
        \tileAA{1}{0}{1}{90}{black}{-15}{0.8}{2};
        \tileAA{1}{1}{1}{0}{black}{15}{0.8}{2};
        \tileAA{1}{2}{1}{0}{black}{-15}{0.8}{2};
        \tileAA{1}{3}{1}{90}{black}{15}{0.8}{2};
        \tileAA{2}{0}{1}{0}{black}{15}{0.8}{2};
        \tileAA{2}{1}{1}{90}{black}{-15}{0.8}{2};
        \tileAA{2}{2}{1}{90}{black}{15}{0.8}{2};
        \tileAA{2}{3}{1}{0}{black}{-15}{0.8}{2};
\end{scope}

\begin{scope}[shift = {(8,0)}]
        \tileAA{0}{0}{1}{0}{black}{45}{0.71}{2};
        \tileAA{0}{1}{1}{0}{black}{-45}{0.71}{2};
        \tileAA{0}{2}{1}{0}{black}{45}{0.71}{2};
        \tileAA{0}{3}{1}{90}{black}{-45}{0.71}{2};
        \tileAA{1}{0}{1}{90}{black}{-45}{0.71}{2};
        \tileAA{1}{1}{1}{0}{black}{45}{0.71}{2};
        \tileAA{1}{2}{1}{0}{black}{-45}{0.71}{2};
        \tileAA{1}{3}{1}{90}{black}{45}{0.71}{2};
        \tileAA{2}{0}{1}{0}{black}{45}{0.71}{2};
        \tileAA{2}{1}{1}{90}{black}{-45}{0.71}{2};
        \tileAA{2}{2}{1}{90}{black}{45}{0.71}{2};
        \tileAA{2}{3}{1}{0}{black}{-45}{0.71}{2};
\end{scope}

\begin{scope}[shift = {(12,0)}]
        \tileAA{0}{0}{1}{0}{black}{45}{0.71}{2};
        \tileAA{0}{1}{1}{0}{black}{-45}{0.71}{2};
        \tileAA{0}{2}{1}{0}{black}{45}{0.71}{2};
        \tileAA{0}{3}{1}{90}{black}{-45}{0.71}{2};
        \tileAA{1}{0}{1}{90}{black}{-45}{0.71}{2};
        \tileAA{1}{1}{1}{0}{black}{45}{0.71}{2};
        \tileAA{1}{2}{1}{0}{black}{-45}{0.71}{2};
        \tileAA{1}{3}{1}{90}{black}{45}{0.71}{2};
        \tileAA{2}{0}{1}{0}{black}{45}{0.71}{2};
        \tileAA{2}{1}{1}{90}{black}{-45}{0.71}{2};
        \tileAA{2}{2}{1}{90}{black}{45}{0.71}{2};
        \tileAA{2}{3}{1}{0}{black}{-45}{0.71}{2};

        \tileAAA{0.5}{0.5}{1}{0}{black}{-45}{0.71}{2};
        \tileAAA{0.5}{1.5}{1}{90}{black}{-45}{0.71}{2};
        \tileAAA{0.5}{2.5}{1}{0}{black}{-45}{0.71}{2};

        \tileAAA{1.5}{0.5}{1}{0}{black}{45}{0.71}{2};
        \tileAAA{1.5}{1.5}{1}{90}{black}{45}{0.71}{2};
        \tileAAA{1.5}{2.5}{1}{0}{black}{45}{0.71}{2};                                
\end{scope}

\draw[->](3.25,2)--(3.75,2);
\draw[->](7.25,2)--(7.75,2);
\draw[->](11.25,2)--(11.75,2);

\node at (16,2){OR};

\begin{scope}[shift={(17,0)}]
        \tileAA{0}{0}{1}{0}{black}{0}{1}{4};
        \tileAA{0}{1}{1}{0}{black}{0}{1}{4};
        \tileAA{0}{2}{1}{0}{black}{0}{1}{4};
        \tileAA{0}{3}{1}{90}{black}{0}{1}{4};
        \tileAA{1}{0}{1}{90}{black}{0}{1}{4};
        \tileAA{1}{1}{1}{0}{black}{0}{1}{4};
        \tileAA{1}{2}{1}{0}{black}{0}{1}{4};
        \tileAA{1}{3}{1}{90}{black}{0}{1}{4};
        \tileAA{2}{0}{1}{0}{black}{0}{1}{4};
        \tileAA{2}{1}{1}{90}{black}{0}{1}{4};
        \tileAA{2}{2}{1}{90}{black}{0}{1}{4};
        \tileAA{2}{3}{1}{0}{black}{0}{1}{4};

\begin{scope}[shift = {(4,0)}]
        \tileAA{0}{0}{1}{0}{black}{-15}{0.8}{2};
        \tileAA{0}{1}{1}{0}{black}{15}{0.8}{2};
        \tileAA{0}{2}{1}{0}{black}{-15}{0.8}{2};
        \tileAA{0}{3}{1}{90}{black}{15}{0.8}{2};
        \tileAA{1}{0}{1}{90}{black}{15}{0.8}{2};
        \tileAA{1}{1}{1}{0}{black}{-15}{0.8}{2};
        \tileAA{1}{2}{1}{0}{black}{15}{0.8}{2};
        \tileAA{1}{3}{1}{90}{black}{-15}{0.8}{2};
        \tileAA{2}{0}{1}{0}{black}{-15}{0.8}{2};
        \tileAA{2}{1}{1}{90}{black}{15}{0.8}{2};
        \tileAA{2}{2}{1}{90}{black}{-15}{0.8}{2};
        \tileAA{2}{3}{1}{0}{black}{15}{0.8}{2};
\end{scope}

\begin{scope}[shift = {(8,0)}]
        \tileAA{0}{0}{1}{0}{black}{-45}{0.71}{2};
        \tileAA{0}{1}{1}{0}{black}{45}{0.71}{2};
        \tileAA{0}{2}{1}{0}{black}{-45}{0.71}{2};
        \tileAA{0}{3}{1}{90}{black}{45}{0.71}{2};
        \tileAA{1}{0}{1}{90}{black}{45}{0.71}{2};
        \tileAA{1}{1}{1}{0}{black}{-45}{0.71}{2};
        \tileAA{1}{2}{1}{0}{black}{45}{0.71}{2};
        \tileAA{1}{3}{1}{90}{black}{-45}{0.71}{2};
        \tileAA{2}{0}{1}{0}{black}{-45}{0.71}{2};
        \tileAA{2}{1}{1}{90}{black}{45}{0.71}{2};
        \tileAA{2}{2}{1}{90}{black}{-45}{0.71}{2};
        \tileAA{2}{3}{1}{0}{black}{45}{0.71}{2};
\end{scope}

\begin{scope}[shift = {(12,0)}]
        \tileAA{0}{0}{1}{0}{black}{-45}{0.71}{2};
        \tileAA{0}{1}{1}{0}{black}{45}{0.71}{2};
        \tileAA{0}{2}{1}{0}{black}{-45}{0.71}{2};
        \tileAA{0}{3}{1}{90}{black}{45}{0.71}{2};
        \tileAA{1}{0}{1}{90}{black}{45}{0.71}{2};
        \tileAA{1}{1}{1}{0}{black}{-45}{0.71}{2};
        \tileAA{1}{2}{1}{0}{black}{45}{0.71}{2};
        \tileAA{1}{3}{1}{90}{black}{-45}{0.71}{2};
        \tileAA{2}{0}{1}{0}{black}{-45}{0.71}{2};
        \tileAA{2}{1}{1}{90}{black}{45}{0.71}{2};
        \tileAA{2}{2}{1}{90}{black}{-45}{0.71}{2};
        \tileAA{2}{3}{1}{0}{black}{45}{0.71}{2};

        \tileAAA{0.5}{0.5}{1}{90}{black}{-45}{0.71}{2};
        \tileAAA{0.5}{1.5}{1}{0}{black}{-45}{0.71}{2};
        \tileAAA{0.5}{2.5}{1}{90}{black}{-45}{0.71}{2};

        \tileAAA{1.5}{0.5}{1}{90}{black}{45}{0.71}{2};
        \tileAAA{1.5}{1.5}{1}{0}{black}{45}{0.71}{2};
        \tileAAA{1.5}{2.5}{1}{90}{black}{45}{0.71}{2};                                
\end{scope}

\draw[->](3.25,2)--(3.75,2);
\draw[->](7.25,2)--(7.75,2);
\draw[->](11.25,2)--(11.75,2);

  \end{scope}

      \end{tikzpicture}
      \caption{Hinged Truchet tiling, two possible directions}
        \label{hinge}
      \end{figure}

    \section{Hinged Truchet curve generation}
    
    Figure~\ref{hinge} shows how a hinged tiling application to a Truchet tiling leads
    to another Truchet tiling.
    After the hinging process is applied to open the tiling to maximum extent
    (third picture in sequences), there is only
    one way in which the tiling can be completed with the same kind of tiles to preserve
    the connected components of the tiling, as shown.  The added tiles
    are coloured grey instead of black to help distinguish them.
    There are two possible ways to hinge the tiling, depending on the direction
    of rotation.  Alternate squares rotate in opposite directions, and
    remain connected at appropriate corners.  In a traditional hinged tiling,
    the squares would remain the same size, but in this version, each square
    tile is scaled down so that the vertices which meet
    adjacent tiles move along  the original tile.
    Figure~\ref{additions} left
    shows how the added tiles are inserted, in either
    case of the direction of turning.  The tiles are coloured according to
    the direction of rotation.
    These two cases are determined by the convention that the tiles are
    labeled alternately odd or even, which can be considered as the
    parity of the sum of the coordinates when tiles are indexed
    by integer coordinates.  Operation 0 is where even tiles rotate
    clockwise, odd tiles anti-clockwise; operation 1 is the reverse.
    After the operation is applied, I use the convention that all the
    old tiles are now even, and all the added tiles are now odd.
    It is important to have some consistent convention, in order to
    apply the operations repeatedly.

\newcommand\twistC[1]{
         \clip (0,0) rectangle (1,1);
          \begin{scope}[shift = {(0.5,0.5)}]
          \begin{scope}[scale = {0.7071},rotate={45}]
          \begin{scope}[shift = {(-0.5,-0.5)}]
            \tileAA{1}{0}{1}{90}{gray!50!white}{0}{1}{4};
            \tileAA{0}{1}{1}{0}{gray!50!white}{0}{1}{4};
            \tileAA{0}{-1}{1}{0}{gray!50!white}{0}{1}{4};
            \tileAA{-1}{0}{1}{90}{gray!50!white}{0}{1}{4};
            \filldraw[cyan!50!white](0,0) rectangle (1,1);
            \draw[->] (0.75,0.5) arc (0:180:0.25);            
         \end{scope}
          \end{scope}
          \end{scope}
            \draw[gray](0,0) rectangle (1,1);
}

\newcommand\twistD[1]{
         \clip (0,0) rectangle (1,1);
          \begin{scope}[shift = {(0.5,0.5)}]
          \begin{scope}[scale = {0.7071},rotate={45}]
          \begin{scope}[shift = {(-0.5,-0.5)}]
            \tileAA{1}{0}{1}{0}{gray!50!white}{0}{1}{4};
            \tileAA{0}{1}{1}{90}{gray!50!white}{0}{1}{4};
            \tileAA{0}{-1}{1}{90}{gray!50!white}{0}{1}{4};
            \tileAA{-1}{0}{1}{0}{gray!50!white}{0}{1}{4};
            \filldraw[orange!50!white](0,0) rectangle (1,1);
            \draw[->] (0.25,0.5) arc (180:0:0.25);
         \end{scope}
          \end{scope}
          \end{scope}
          \draw[gray](0,0) rectangle (1,1);          
          }

    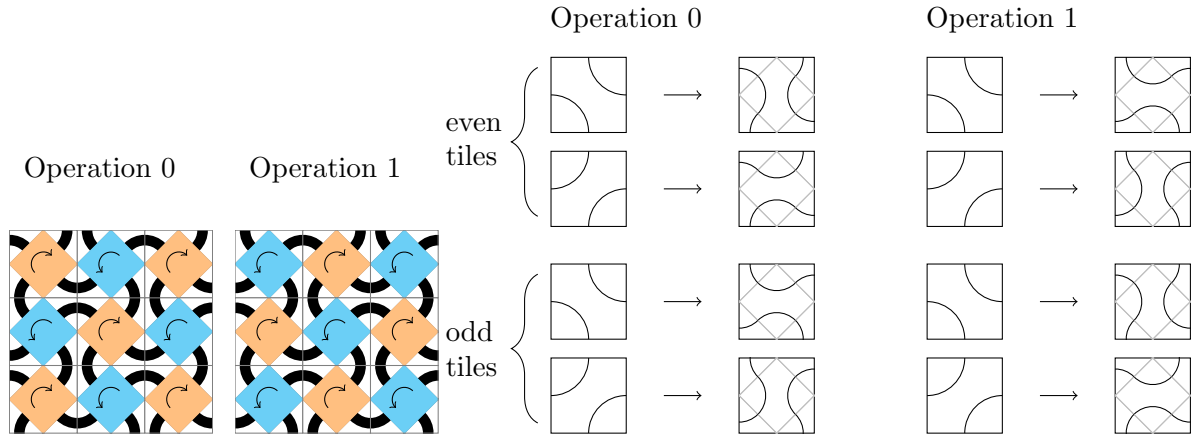
\begin{figure}
      \begin{tikzpicture}
        \begin{scope}[shift={(-5.7,0)},scale = {0.9}]
          \begin{scope}[shift={(0,1)}]\twistC{0}\end{scope}
          \begin{scope}[shift={(1,0)}]\twistC{0}\end{scope}
          \begin{scope}[shift={(1,2)}]\twistC{0}\end{scope}
          \begin{scope}[shift={(2,1)}]\twistC{0}\end{scope}          
          \begin{scope}[shift={(1,1)}]\twistD{0}\end{scope}
          \begin{scope}[shift={(2,0)}]\twistD{0}\end{scope}
          \begin{scope}[shift={(0,0)}]\twistD{0}\end{scope}
          \begin{scope}[shift={(0,2)}]\twistD{0}\end{scope}
          \begin{scope}[shift={(2,2)}]\twistD{0}\end{scope}
        \end{scope}

        \begin{scope}[shift={(-2.7,0)},scale = {0.9}]
          \begin{scope}[shift={(0,1)}]\twistD{0}\end{scope}
          \begin{scope}[shift={(1,0)}]\twistD{0}\end{scope}
          \begin{scope}[shift={(1,2)}]\twistD{0}\end{scope}
          \begin{scope}[shift={(2,1)}]\twistD{0}\end{scope}          
          \begin{scope}[shift={(1,1)}]\twistC{0}\end{scope}
          \begin{scope}[shift={(2,0)}]\twistC{0}\end{scope}
          \begin{scope}[shift={(0,0)}]\twistC{0}\end{scope}
          \begin{scope}[shift={(0,2)}]\twistC{0}\end{scope}
          \begin{scope}[shift={(2,2)}]\twistC{0}\end{scope}
          \end{scope}

        \node at(-4.5,3.5){Operation 0};
        \node at(-1.5,3.5){Operation 1};
        \node at(2.5,5.5){Operation 0};
        \node at(7.5,5.5){Operation 1};        
        \tileA{1.5}{4}{1}{90}{black};
        \tileA{1.5}{2.75}{1}{0}{black};
        \tileA{1.5}{1.25}{1}{90}{black};
        \tileA{1.5}{0}{1}{0}{black};
        \node (a) at (1.5,0){};
        \node (b) at (1.5,2.25){};
        \node (c) at (1.5,2.75){};
        \node (d) at (1.5,5){};        
\draw  [decorate,decoration={brace,amplitude=10pt,raise=5pt}]
(a) -- (b) node [black,midway,yshift=0.cm,xshift=-.4cm]
       {\parbox{2cm}{odd\\ tiles}};

\draw  [decorate,decoration={brace,amplitude=10pt,raise=5pt}]
(c) -- (d) node [black,midway,yshift=0.cm,xshift=-.4cm]
       {\parbox{2cm}{even\\ tiles}};

       \draw[->](3,0.5)--(3.5,0.5);              
       \draw[->](3,1.75)--(3.5,1.75);              
       \draw[->](3,3.25)--(3.5,3.25);              
\draw[->](3,4.5)--(3.5,4.5);
       
\begin{scope}[shift = {(4,0)}]\twistA{0}\end{scope}
\begin{scope}[shift = {(4,1.25)}]\twistA{90}\end{scope}
\begin{scope}[shift = {(4,2.75)}]\twistB{90}\end{scope}
\begin{scope}[shift = {(4,4.)}]\twistB{0}\end{scope}

\begin{scope}[shift = {(5,0)}]
\tileA{1.5}{4}{1}{90}{black};
        \tileA{1.5}{2.75}{1}{0}{black};
        \tileA{1.5}{1.25}{1}{90}{black};
        \tileA{1.5}{0}{1}{0}{black};

       \draw[->](3,0.5)--(3.5,0.5);              
       \draw[->](3,1.75)--(3.5,1.75);              
       \draw[->](3,3.25)--(3.5,3.25);              
\draw[->](3,4.5)--(3.5,4.5);
       
\begin{scope}[shift = {(4,0)}]\twistB{90}\end{scope}
\begin{scope}[shift = {(4,1.25)}]\twistB{0}\end{scope}
\begin{scope}[shift = {(4,2.75)}]\twistA{0}\end{scope}
\begin{scope}[shift = {(4,4.)}]\twistA{90}\end{scope}

        \end{scope}
       
        \end{tikzpicture}
      \caption{How the new tiles are added in hinged tiling, for either of the two operations, and an alternative
        way to see the transformations as a replacement rule.  Even and odd
      refer to the grid position of the tile.}
      \label{additions}
      \end{figure}
    
    This procedure can also be viewed as a replacement of
    each tile by smaller tiles as
    in Figure~\ref{additions} right.
    This can be compared with the replacement rules for the Hilbert
    space filling curve, explored in detail in \cite{Ozkaraca},
    however, in this case the tile is divided into parts which do not
    fit as single  whole units into the original tile.  Possibly there are
    other dissections possible for tile replacements.
    The replacement rule in Figure~\ref{additions} only has four cases
    because the tiles have $180^\circ$ rotational symmetry.  A longer
    list of rules would be needed in less symmetric cases.
    The tiles can be rotated in a continuous fashion, and the added tiles
    can be incorporated as the original tiles are rotated, rather than
    added at the end,    so the 
    procedure can be considered as a
    continuous transformation rather than a discrete operation.
    This shows that the number of components of a Truchet tiling path
    configuration does not
    change with this operation.
    A program illustrating this in motion can be found at \cite{Verrill}.

    \section{L-systems}

    The effect of the procedure described on the curves
    can be described as an L-system.
    An L-system is a method of describing growth processes
    developed by Lindenmayer
    \cite{LINDENMAYER1968300}.
    In our case, words describe directions for
    drawing a path.
    A replacement rule on the letters of the word gives instructions
    to generate the next step in the procedure.

    Paths on Truchet tilings of the form we are interested in can be given
    consistent ``in'' and ``out'' directions
    as shown by the arrows in the example in
    Figure~\ref{directions}.  The tiles' orientation determine whether the
    path curves left or right into the next square.  Thus each path can
    be described by a sequence such as $RvLhRvLhLvR$ as seen in the figure,
    and more generally, a path is given by a sequence $XvXhXvXh...$ where
    $X$ is either $R$ or $L$.   The $R$ label is used when the path is turning
    right, and the $L$ when turning left.
    In the figure, for less clutter and more clarity, the letters $R$ and $L$
    only appear once per tile, which is possible since both paths in any tile
    turn in the same direction.
    The label $v$ is used for crossing
    a vertical line, and $h$ for crossing a horizontal line.  In the figure,
    the dark squares are considered even, and this is an application of operation 0.  As the operation is applied, the new tiles added in the spaces labeled
    $v$ or $h$ are replaced with tiles with paths turning either left or right
    respectively.  When the new tiles appear, there will be new ``horizontal'' and
    ``vertical'' lines to cross.  The tiling is now rotated by $45^\circ$.
    However, we are using the convention that all the original tiles remain or
    become even, and so this determines which of the diagonal lines are
    labeled $h$ or $v$, because the direction of flow is always into an even
    tile from a $h$ line, and into an odd tile from a $v$ line.  So in the
    rightmost diagram in Figure~\ref{directions} the NW--SE lines should be
    labeled $h$ and the NE--SW lines should be labeled $v$.

            \newcommand\tileCCC[8]{
        \begin{scope}[shift = {(#1+#3/2,#2+#3/2)}]
\begin{scope}[scale={#7}]
          \begin{scope}[rotate={#6}]
          \draw[gray, fill=#5](-#3/2,-#3/2) rectangle (#3/2,#3/2);
    \begin{scope}[rotate = {#4}]
      \draw[<-,line width = #8] (-#3/2,0) arc (-90:0:#3/2);
\draw[->,line width = #8] (0,-#3/2) arc (180:90:#3/2);    
    \end{scope}
    \end{scope}    
        \end{scope}
      \end{scope}
            }

            \newcommand\tileDDD[8]{
        \begin{scope}[shift = {(#1+#3/2,#2+#3/2)}]
\begin{scope}[scale={#7}]
          \begin{scope}[rotate={#6}]
          \draw[gray, fill=#5](-#3/2,-#3/2) rectangle (#3/2,#3/2);
    \begin{scope}[rotate = {#4}]
      \draw[->,line width = #8] (-#3/2,0) arc (-90:0:#3/2);
\draw[<-,line width = #8] (0,-#3/2) arc (180:90:#3/2);    
    \end{scope}
    \end{scope}    
        \end{scope}
      \end{scope}
  }              

            \newcommand\tileAB[5]{
  \begin{scope}[shift = {(#1+#3/2,#2+#3/2)}]
    \begin{scope}[rotate = {#4}]
    \draw[gray,line width = 2] (-#3/2,0) arc (-90:0:#3/2);
\draw[gray,line width = 2] (0,-#3/2) arc (180:90:#3/2);    
    \end{scope}
      \end{scope}
  }

\newcommand\twistBB[1]{
         \clip (0,0) rectangle (1,1);
          \begin{scope}[shift = {(0.5,0.5)}]
          \begin{scope}[scale = {0.86},rotate={-10}]
          \begin{scope}[shift = {(-0.5,-0.5)}]
            \tileAB{1}{0}{1}{0}{gray!50!white};
            \tileAB{0}{1}{1}{90}{gray!50!white};
            \tileAB{0}{-1}{1}{90}{gray!50!white};
            \tileAB{-1}{0}{1}{0}{gray!50!white};            
         \end{scope}
          \end{scope}
                  \end{scope}
          }

\newcommand\twistAA[1]{
         \clip (0,0) rectangle (1,1);
          \begin{scope}[shift = {(0.5,0.5)}]
          \begin{scope}[scale = {0.86},rotate={10}]
          \begin{scope}[shift = {(-0.5,-0.5)}]
            \tileAB{0}{0}{1}{#1}{gray!50!white};
            \tileAB{1}{0}{1}{90}{gray!50!white};
            \tileAB{0}{1}{1}{0}{gray!50!white};
            \tileAB{0}{-1}{1}{0}{gray!50!white};
            \tileAB{-1}{0}{1}{90}{gray!50!white};            
         \end{scope}
          \end{scope}
          \end{scope}          
}

\newcommand\twistBC[1]{
         \clip (0,0) rectangle (1,1);
          \begin{scope}[shift = {(0.5,0.5)}]
          \begin{scope}[scale = {0.71},rotate={-45}]
          \begin{scope}[shift = {(-0.5,-0.5)}]
            \tileAB{1}{0}{1}{0}{gray!50!white};
            \tileAB{0}{1}{1}{90}{gray!50!white};
            \tileAB{0}{-1}{1}{90}{gray!50!white};
            \tileAB{-1}{0}{1}{0}{gray!50!white};            
         \end{scope}
          \end{scope}
                  \end{scope}
          }

\newcommand\twistAC[1]{
         \clip (0,0) rectangle (1,1);
          \begin{scope}[shift = {(0.5,0.5)}]
          \begin{scope}[scale = {0.71},rotate={45}]
          \begin{scope}[shift = {(-0.5,-0.5)}]
            \tileAB{0}{0}{1}{#1}{gray!50!white};
            \tileAB{1}{0}{1}{90}{gray!50!white};
            \tileAB{0}{1}{1}{0}{gray!50!white};
            \tileAB{0}{-1}{1}{0}{gray!50!white};
            \tileAB{-1}{0}{1}{90}{gray!50!white};            
         \end{scope}
          \end{scope}
          \end{scope}          
}

        \begin{figure}
          \begin{tikzpicture}[scale={1.4}]

            \tileCCC{0}{0}{1}{0}{pink!50!white}{0}{1}{2}
            \tileCCC{2}{0}{1}{0}{pink!50!white}{0}{1}{2}            
            \tileDDD{0}{2}{1}{90}{pink!50!white}{0}{1}{2}
            \tileDDD{1}{1}{1}{90}{pink!50!white}{0}{1}{2}
            \tileDDD{2}{2}{1}{90}{pink!50!white}{0}{1}{2}

            \tileDDD{2}{1}{1}{90}{yellow!60!white}{90}{1}{2}
            \tileDDD{1}{2}{1}{90}{yellow!60!white}{90}{1}{2}
            \tileDDD{0}{1}{1}{90}{yellow!60!white}{90}{1}{2}
            \tileCCC{1}{0}{1}{0}{yellow!60!white}{90}{1}{2}

            \node at (0.5,0.5){$R$};
            \node at (0.5,1.5){$L$};
            \node at (0.5,2.5){$L$};
            \node at (1.5,0.5){$R$};
            \node at (1.5,1.5){$L$};
            \node at (1.5,2.5){$L$};
            \node at (2.5,0.5){$R$};
            \node at (2.5,1.5){$L$};
            \node at (2.5,2.5){$L$};

            \node at (0,0.4){$_v$}; \node at (0.4,0){$_h$};
            \node at (1,0.4){$_v$}; \node at (0.4,1){$_h$};
            \node at (2,0.4){$_v$}; \node at (0.4,2){$_h$};
            \node at (3,0.4){$_v$}; \node at (0.4,3){$_h$};
            \node at (3,0.4){$_v$};  \node at (3,0.4){$_h$};            
            \node at (0,1.4){$_v$}; \node at (1.4,0){$_h$};
            \node at (1,1.4){$_v$}; \node at (1.4,1){$_h$};
            \node at (2,1.4){$_v$}; \node at (1.4,2){$_h$};
            \node at (3,1.4){$_v$}; \node at (1.4,3){$_h$};
            \node at (3,1.4){$_v$}; \node at (3,1.4){$_h$};            
            \node at (0,2.4){$_v$}; \node at (2.4,0){$_h$};
            \node at (1,2.4){$_v$}; \node at (2.4,1){$_h$};
            \node at (2,2.4){$_v$}; \node at (2.4,2){$_h$};
            \node at (3,2.4){$_v$}; \node at (2.4,3){$_h$};

            \begin{scope}[shift={(4,0)}]
              \begin{scope}
              \clip (0,0) rectangle (3,3);

              \begin{scope}[shift={(0,0)}]\twistBB{0}\end{scope}
              \begin{scope}[shift={(1,1)}]\twistBB{0}\end{scope}              
              \begin{scope}[shift={(2,2)}]\twistBB{0}\end{scope}
              \begin{scope}[shift={(2,0)}]\twistBB{0}\end{scope}              
              \begin{scope}[shift={(0,2)}]\twistBB{0}\end{scope}
              \begin{scope}[shift={(0,1)}]\twistAA{0}\end{scope}
              \begin{scope}[shift={(1,0)}]\twistAA{0}\end{scope}
              \begin{scope}[shift={(2,1)}]\twistAA{0}\end{scope}
              \begin{scope}[shift={(1,2)}]\twistAA{0}\end{scope}

            \tileCCC{0}{0}{1}{0}{pink!50!white}{-10}{0.86}{2}
            \tileCCC{2}{0}{1}{0}{pink!50!white}{-10}{0.86}{2}            
            \tileDDD{0}{2}{1}{90}{pink!50!white}{-10}{0.86}{2}
            \tileDDD{1}{1}{1}{90}{pink!50!white}{-10}{0.86}{2}
            \tileDDD{2}{2}{1}{90}{pink!50!white}{-10}{0.86}{2}

            \tileDDD{2}{1}{1}{0}{yellow!60!white}{10}{0.86}{2}
            \tileDDD{1}{2}{1}{0}{yellow!60!white}{10}{0.86}{2}
            \tileDDD{0}{1}{1}{0}{yellow!60!white}{10}{0.86}{2}
            \tileCCC{1}{0}{1}{90}{yellow!60!white}{10}{0.86}{2}

            \node at (0.5,0.5){$R$};
            \node at (0.5,1.5){$L$};
            \node at (0.5,2.5){$L$};
            \node at (1.5,0.5){$R$};
            \node at (1.5,1.5){$L$};
            \node at (1.5,2.5){$L$};
            \node at (2.5,0.5){$R$};
            \node at (2.5,1.5){$L$};
            \node at (2.5,2.5){$L$};            
            \end{scope}
              \node at(0,1){$v$};
              \node at(0,2){$h$};
              \node at(1,0){$v$};
              \node at(1,1){$h$};
              \node at(1,2){$v$};
              \node at(2,0){$h$};
              \node at(2,1){$v$};
              \node at(2,2){$h$};
              \node at(1,3){$h$};                            
              \node at(2,3){$v$};
              \node at(3,1){$h$};
              \node at(3,2){$v$};

            \end{scope}

            \begin{scope}[shift={(8,0)}]
              \begin{scope}
              \clip (0,0) rectangle (3,3);

              \begin{scope}[shift={(0,0)}]\twistBC{0}\end{scope}
              \begin{scope}[shift={(1,1)}]\twistBC{0}\end{scope}              
              \begin{scope}[shift={(2,2)}]\twistBC{0}\end{scope}
              \begin{scope}[shift={(2,0)}]\twistBC{0}\end{scope}              
              \begin{scope}[shift={(0,2)}]\twistBC{0}\end{scope}
              \begin{scope}[shift={(0,1)}]\twistAC{0}\end{scope}
              \begin{scope}[shift={(1,0)}]\twistAC{0}\end{scope}
              \begin{scope}[shift={(2,1)}]\twistAC{0}\end{scope}
              \begin{scope}[shift={(1,2)}]\twistAC{0}\end{scope}

            \tileCCC{0}{0}{1}{0}{pink!50!white}{-45}{0.71}{2}
            \tileCCC{2}{0}{1}{0}{pink!50!white}{-45}{0.71}{2}            
            \tileDDD{0}{2}{1}{90}{pink!50!white}{-45}{0.71}{2}
            \tileDDD{1}{1}{1}{90}{pink!50!white}{-45}{0.71}{2}
            \tileDDD{2}{2}{1}{90}{pink!50!white}{-45}{0.71}{2}

            \tileDDD{2}{1}{1}{0}{yellow!60!white}{45}{0.71}{2}
            \tileDDD{1}{2}{1}{0}{yellow!60!white}{45}{0.71}{2}
            \tileDDD{0}{1}{1}{0}{yellow!60!white}{45}{0.71}{2}
            \tileCCC{1}{0}{1}{90}{yellow!60!white}{45}{0.71}{2}

            \node at (0.5,0.5){$R$};
            \node at (0.5,1.5){$L$};
            \node at (0.5,2.5){$L$};
            \node at (1.5,0.5){$R$};
            \node at (1.5,1.5){$L$};
            \node at (1.5,2.5){$L$};
            \node at (2.5,0.5){$R$};
            \node at (2.5,1.5){$L$};
            \node at (2.5,2.5){$L$};            
              \end{scope}
                            \node at(0,1){$L$};
              \node at(0,2){$R$};
              \node at(1,0){$L$};
              \node at(1,1){$R$};
              \node at(1,2){$L$};
              \node at(2,0){$R$};
              \node at(2,1){$L$};
              \node at(2,2){$R$};
              \node at(1,3){$R$};                            
              \node at(2,3){$L$};
              \node at(3,1){$R$};
              \node at(3,2){$L$};      
            
\end{scope}
            
          \end{tikzpicture}
      \caption{Sequence showing transformation of paths,
        directions of flow, and labeling of paths on a Truchet tiling.
      This shows operation $0$, with $v\rightarrow L$ and $h\rightarrow R$.}
      \label{directions}
    \end{figure}
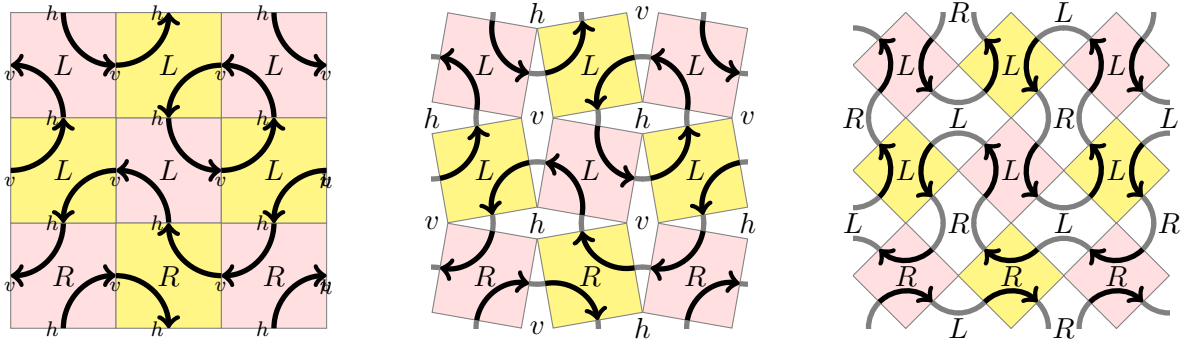

    The rules for the procedures described in this paper are as in
    Table~\ref{tab:L}

      \begin{table}
     \begin{tabular}{lll|lll}
       \multicolumn{3}{c}{operation 0} &        \multicolumn{3}{c}{operation 1}\\
         L & $\rightarrow$ &          L & L & $\rightarrow$ & L\\
         R & $\rightarrow$ &          R & R & $\rightarrow$ & R\\
         h & $\rightarrow$ &          hLv & h & $\rightarrow$ & hRv\\
         v & $\rightarrow$ &          hRv & v &$\rightarrow$ & hLv\\
     \end{tabular}
     \caption{L-system description of hinged tiling operations.}
     \label{tab:L}
      \end{table}

      \subsection{Notation}
      For the remainder of this paper, I will denote paths by a sequence
      of $L$s and $R$s, e.g., $LRRRRLL$, and combinations of operations $0$
      and $1$ by a sequence of $0$s and $1$s.  Thus for example
      $00101LLR$ means start with a path defined by the directions
      $LLR$, then apply the operations in the order from {\it{left to
        right}}, i.e., operation $0$, then operation $0$, then operation $1$ and so on.
    
    \section{Comparison with Heighway's Fractal dragon}
\begin{figure}
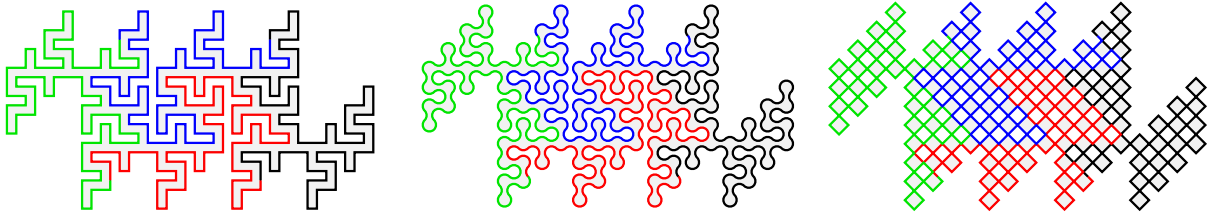

    \begin{tikzpicture}      
      \input{0100010lllldiagonal.tex}
\begin{scope}[shift={(9.8,3.8)}]
      \input{0100010llllcurves.tex}      
\end{scope}
\begin{scope}[shift={(11,0.05)}]
      \input{0100010llllangle.tex}      
\end{scope}
    \end{tikzpicture}
    \caption{Comparison of $0100010RRRR$ for different tile designs.}
    \label{compare}
\end{figure}

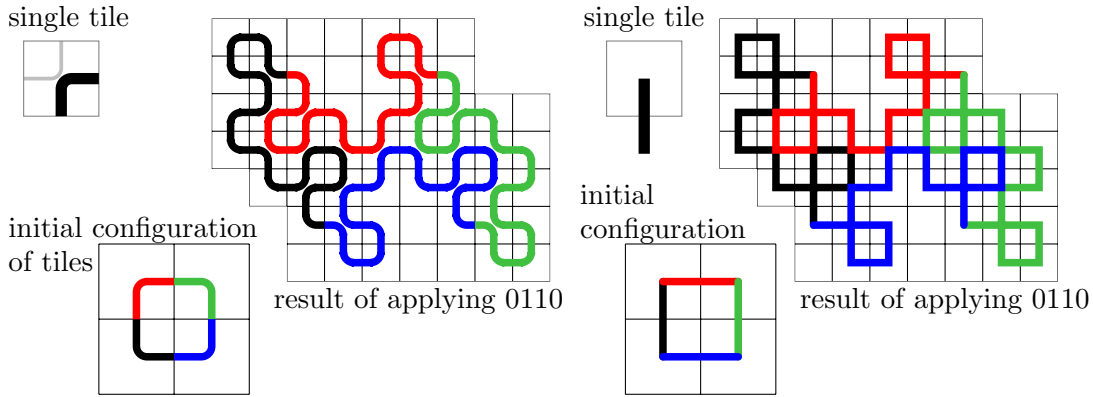
\begin{figure}
  \begin{tikzpicture}
   \tileCC{-2}{2.7}{1}{0}{gray}
    \node at (-1.4,4){single tile};
    \node at (-0.2,1){\parbox{4cm}{initial configuration\\ of tiles}};
    \node at (7.4,1.4){\parbox{4cm}{initial\\ configuration}};   
    \node at (3.25,0.25){result of applying $0110$};
    \node at (10.25,0.25){result of applying $0110$};            
    \node at (6.25,4){single tile};    
    \tileD{5.75}{2.7}{1}{0}{gray}
    \draw (-1,-1) grid (1,1);
    \draw (6,-1) grid (8,1);    
\foreach \i in {0,6.75}{
    \begin{scope}[shift = {(\i,0)}]
      \clip (0.5,2)--(0.5,4)--(4,4)--(4,3)--(5,3)--(5,0.5)--(1.5,0.5)--(1.5,1.5)
      --(1,1.5)--(1,2) -- cycle;
      \draw[step=.5] (0.5,0.5) grid (5,4);
    \end{scope}
}
    \begin{scope}[shift={(0,0.5)}]    
      \begin{scope}[scale={0.005}]
\begin{scope}[rotate={180}]
\begin{scope}[shift={(0,200)}]
\draw[rounded corners,black,  line width=1mm] (0,0) -- ++(100,0)-- ++(0,-100.000000) ;
\draw[rounded corners,red,  line width=1mm] (100,-100) -- ++(0,-100.000000)-- ++(-100,0);
\draw[rounded corners,green!50!gray,  line width=1mm] (0,-200) -- ++(-100,0)-- ++(0,100.000000) ;
\draw[rounded corners,blue,  line width=1mm] (-100,-100) -- ++(0,100.000000)-- ++(100,0);
\end{scope}
\end{scope}
\end{scope}
      \end{scope}
    \begin{scope}[shift={(7.5,0.5)}]    
      \begin{scope}[scale={0.005}]
\begin{scope}[rotate={180}]
\begin{scope}[shift={(200,200)}]
\draw[black, line width=1mm,line cap=round] (0,0) -- ++(0,-200.000000) ;
\draw[red, line width=1mm,line cap=round] (0,-200) -- ++(-200,0.000000) ;
\draw[green!50!gray, line width=1mm,line cap=round] (-200,-200) -- ++(0,200.000000) ;
\draw[blue, line width=1mm,line cap=round] (-200,0) -- ++(200,0.000000) ;
\end{scope}
\end{scope}
\end{scope}
\end{scope}
    \begin{scope}[shift={(3.5,3.25)}]
  \begin{scope}[scale={0.01}]
\begin{scope}[rotate={180}]
\begin{scope}[shift={(149.99999999999991,199.99999999999994)}]
\draw[rounded corners, black,  line width=1mm] (0,0) -- ++(25,0)-- ++(0,-25.000000) -- ++(0,-25.000000)-- ++(-25,0)-- ++(-25,0)-- ++(0,-25.000000) -- ++(0,-25.000000)-- ++(25,0)-- ++(25,0)-- ++(0,25.000000) -- ++(0,25.000000)-- ++(25,0)-- ++(25,0)-- ++(0,-25.000000) -- ++(0,-25.000000)-- ++(25,0)-- ++(25,0)-- ++(0,-25.000000) -- ++(0,-25.000000)-- ++(-25,0)-- ++(-25,0)-- ++(0,-25.000000) -- ++(0,-25.000000)-- ++(25,0)-- ++(25,0)-- ++(0,-25.000000) -- ++(0,-25.000000)-- ++(-25,0)-- ++(-25,0)-- ++(0,25.000000) -- ++(0,25.000000)-- ++(-25,0);
\draw[rounded corners,red,  line width=1mm] (49.999999999999986,-199.99999999999994) -- ++(-25,0)-- ++(0,25.000000) -- ++(0,25.000000)-- ++(25,0)-- ++(25,0)-- ++(0,25.000000) -- ++(0,25.000000)-- ++(-25,0)-- ++(-25,0)-- ++(0,-25.000000) -- ++(0,-25.000000)-- ++(-25,0)-- ++(-25,0)-- ++(0,25.000000) -- ++(0,25.000000)-- ++(-25,0)-- ++(-25,0)-- ++(0,-25.000000) -- ++(0,-25.000000)-- ++(-25,0)-- ++(-25,0)-- ++(0,-25.000000) -- ++(0,-25.000000)-- ++(25,0)-- ++(25,0)-- ++(0,-25.000000) -- ++(0,-25.000000)-- ++(-25,0)-- ++(-25,0)-- ++(0,25.000000) -- ++(0,25.000000)-- ++(-25,0);
\draw[rounded corners,green!50!gray,  line width=1mm] (-149.99999999999991,-199.99999999999994) -- ++(-25,0)-- ++(0,25.000000) -- ++(0,25.000000)-- ++(25,0)-- ++(25,0)-- ++(0,25.000000) -- ++(0,25.000000)-- ++(-25,0)-- ++(-25,0)-- ++(0,-25.000000) -- ++(0,-25.000000)-- ++(-25,0)-- ++(-25,0)-- ++(0,25.000000) -- ++(0,25.000000)-- ++(-25,0)-- ++(-25,0)-- ++(0,25.000000) -- ++(0,25.000000)-- ++(25,0)-- ++(25,0)-- ++(0,25.000000) -- ++(0,25.000000)-- ++(-25,0)-- ++(-25,0)-- ++(0,25.000000) -- ++(0,25.000000)-- ++(25,0)-- ++(25,0)-- ++(0,-25.000000) -- ++(0,-25.000000)-- ++(25,0);
\draw[rounded corners,blue,  line width=1mm] (-199.9999999999999,-2.1316282072803006e-14) -- ++(25,0)-- ++(0,-25.000000) -- ++(0,-25.000000)-- ++(-25,0)-- ++(-25,0)-- ++(0,-25.000000) 
-- ++(0,-25.000000)-- ++(25,0)-- ++(25,0)-- ++(0,25.000000) -- ++(0,25.000000)-- ++(25,0)-- ++(25,0)-- ++(0,-25.000000) -- ++(0,-25.000000)-- ++(25,0)-- ++(25,0)-- ++(0,25.000000) -- ++(0,25.000000)-- ++(25,0)-- ++(25,0)-- ++(0,25.000000) -- ++(0,25.000000)-- ++(-25,0)-- ++(-25,0)-- ++(0,25.000000) -- ++(0,25.000000)-- ++(25,0)-- ++(25,0)-- ++(0,-25.000000) -- ++(0,-25.000000)-- ++(25,0);
\end{scope}
\end{scope}
\end{scope}
\end{scope}
\begin{scope}[shift={(10.5,3.25)}]
  \begin{scope}[scale={0.01}]
\begin{scope}[rotate={180}]
\begin{scope}[shift={(199.99999999999991,199.99999999999991)}]
\draw[black,  line width=1mm,line cap=round] (0,0) -- ++(0,-50.000000) -- ++(-50,0.000000) -- ++(0,-50.000000) -- ++(50,0.000000) -- ++(0,50.000000) -- ++(50,0.000000) -- ++(0,-50.000000) -- ++(50,0.000000) -- ++(0,-50.000000) -- ++(-50,0.000000) -- ++(0,-50.000000) -- ++(50,0.000000) -- ++(0,-50.000000) -- ++(-50,0.000000) -- ++(0,50.000000) -- ++(-50,0.000000) ;
\draw[red,  line width=1mm,line cap=round] (0,-199.99999999999991) -- ++(0,50.000000) -- ++(50,0.000000) -- ++(0,50.000000) -- ++(-50,0.000000) -- ++(0,-50.000000) -- ++(-50,0.000000) -- ++(0,50.000000) -- ++(-50,0.000000) -- ++(0,-50.000000) -- ++(-50,0.000000) -- ++(0,-50.000000) -- ++(50,0.000000) -- ++(0,-50.000000) -- ++(-50,0.000000) -- ++(0,50.000000) -- ++(-50,0.000000) ;
\draw[green!50!gray,  line width=1mm,line cap=round] (-199.99999999999991,-199.99999999999991) -- ++(0,50.000000) -- ++(50,0.000000) -- ++(0,50.000000) -- ++(-50,0.000000) -- ++(0,-50.000000) -- ++(-50,0.000000) -- ++(0,50.000000) -- ++(-50,0.000000) -- ++(0,50.000000) -- ++(50,0.000000) -- ++(0,50.000000) -- ++(-50,0.000000) -- ++(0,50.000000) -- ++(50,0.000000) -- ++(0,-50.000000) -- ++(50,0.000000) ;
\draw[blue,  line width=1mm,line cap=round] (-199.99999999999994,-1.4210854715202004e-14) -- ++(0,-50.000000) -- ++(-50,0.000000) -- ++(0,-50.000000) 
-- ++(50,0.000000) -- ++(0,50.000000) -- ++(50,0.000000) -- ++(0,-50.000000) -- ++(50,0.000000) -- ++(0,50.000000) -- ++(50,0.000000) -- ++(0,50.000000) -- ++(-50,0.000000) -- ++(0,50.000000) -- ++(50,0.000000) -- ++(0,-50.000000) -- ++(50,0.000000) ;
\end{scope}
\end{scope}
\end{scope}
  \end{scope}
        \end{tikzpicture}
  \caption{Comparison of tile in the  hinged tile construction (left), and
    a replacement tile comparable to the
    usual Heighway dragon construction (right).  The sequence of operations
  $0110$ has been applied to the initial square arrangement of line segments.}
    \label{replace}
\end{figure}

    In \cite{Tab} the fractal Heighway dragon curve
    is described as being obtained from
    an L system equivalent to that in Table~\ref{tab:L},
    applying only operation 0 repeatedly.  The rules described in
    \cite{Edgar}
    are interpreted as drawing a line, then
    making a turn, and so on.  For the hinged tiling method, we can put
    different designs on our tiles, which for the finite cases may look
    somewhat different, but in the limit will all end up being the same.
    Figure~\ref{troucheexamples} shows some different possible tiles, and
    Figure~\ref{compare} compares the result of applying the operation
    $0100010$ to the initial configuration $LLLL$.  To compare with the
    construction of the Heighway dragon given in \cite[]{Edgar}, we 
    slightly alter
    our curves as in Figure~\ref{replace}.  All this means is that for
    the usual construction, the lines, if drawn on tiles, are
    considered going from centre of tile
    to centre of tile, whereas in the original
    Truchet tiling construction the tiles
    go from edge to edge.  In the Truchet tiling construction we may have
    many components, producing many curves simultaneously, whereas in the
    standard construction we just start with one component.  So in Figure~\ref{replace} I have coloured the components in different colours, including
    leaving some parts of tiles only having half the design shown.
    
By construction, the Truchet tilings we start with are always
non-self intersecting, and tessellate.  Each step of the iteration always
produces another Truchet tiling.  So results about non-self intersection and
tessellation given in
\cite{Edgar} and \cite{Tab} are now obtained automatically.
    
\section{Fractal dimension}
\label{sec:fractaldim}

\subsection{Space filling}
\begin{figure}
  \begin{tikzpicture}[scale = {0.75}]
    \draw[fill=gray!5!white](0,2) rectangle (2,4);
    \node (a) at (2,1.)
          {\parbox{6cm}{{\tiny
                {single tile containing\\ original line segment}}}};

          \draw[fill=gray!20!white](3,2)--(6,5)--(7,4)--(4,1)--cycle;

    \node (a) at (6.5,0.5)
          {\parbox{6cm}{{\tiny
                {$5$ tiles containing line segment\\ after one step}}}};

          \foreach \i in {0,5}{
            \begin{scope}[shift = {(\i,0)}]
          \draw[fill=gray!20!white](3,4)--(4,5)--(7,2)--(6,1) -- cycle;
          \draw[fill=gray!20!white](3,2)--(6,5)--(7,4)--(4,1)--cycle;
                    \draw[fill=gray!5!white](4,3)--(5,4)--(6,3)--(5,2)--cycle;
          \draw[dotted](4,2) rectangle (6,4);
\end{scope}
          }
\begin{scope}[shift= {(10,3)}]
          \foreach \i in {0,90,180,270}{
            \begin{scope}[rotate = {\i}]
              \begin{scope}[scale={0.5}]
             \draw[fill=gray!20!white](1,1) rectangle (3,3);
            \draw[fill=gray!40!white](-1,1)--(5,1)--  (5,3)--(-1,3)--cycle;
            \draw[fill=gray!40!white](1,-1)--(1,5)--  (3,5)--(3,-1)--cycle;
                            
            \draw(1,1)--(5,1)--  (5,3)--(-1,3)--cycle;
            
            \draw[dashed](2,0)--(0,2)--(2,4)--(4,2)-- cycle;
            \draw[fill=gray!20!white](1,1) rectangle (3,3);
\end{scope}
            \end{scope}
            }
\end{scope}            

    \node (a) at (11.5,0)
          {\parbox{6cm}{{\tiny
                {tiles containing line after 2 steps}}}};

            \draw[dotted](9,2) rectangle (11,4);

    \end{tikzpicture}
  \caption{Argument that after infinite application of $0$ and $1$,
    the resulting curve is contained in a finite region.
First two applications of process shown.
  }
  \label{piles}
  \end{figure}
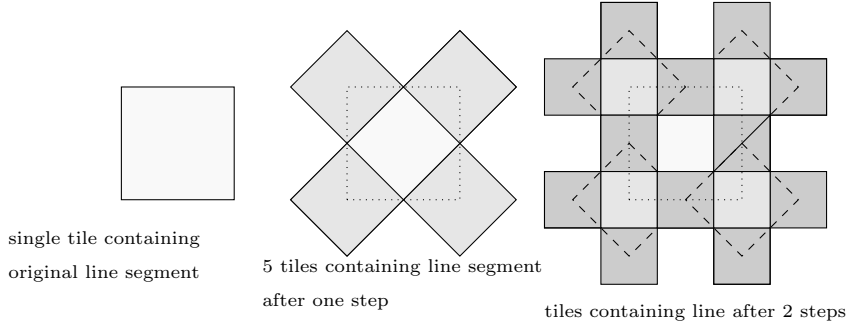

The curves described here are space filling:
take a single tile with just two components.  Consider
the Truchet tiling to take place on the surface of a torus with initially
only one square.  On the torus the configuration now consists of a closed
ring of two components, and no end points.  Now apply any infinite
sequence of the operations $0$ and $1$.  At each step, a finer grid of
tiles is covered.
In the limit the whole torus is covered, with points of
the curve arbitrarily close to any point of the torus.
Now unroll our torus to the plane.  We get a tessellation of the
tiles covering the plane.  We have to rule out the possibility of individual
tiles not being contained in a finite region.
This follows by considering the iteration step in the hinged tiling.
For any single line segment, the iteration step takes the curve to
a new curve completely contained in the same square.  However, each square
is not contained in the previous square, so after one step, a curve is contained
in the original square, and in two adjacent squares,
but the adjacent squares, and the scaled original square have side length only $\sqrt{2}$ the original
side length.  So, taking the initial
squares to have unit length, after two steps,
the curve is contained in a square of
side length $1$ and square of side length $\sqrt{2}$.
The new points added to the curve are at most $\sqrt{2}$ from the initial curve.
At successive steps we can get points further away by increments of $\sqrt{2}$.
So the maximum total possible distance of any point in the limit curve is
$\sum 2^{-n/2}=2-\sqrt{2}<1$.
This is illustrated in Figure~\ref{piles}.
This means that the whole of the initial
square must be covered by fractiles (using the terminology of \cite{Edgar}
for tiles with (possibly) fractal boundary) derived from curves at most in
the adjacent squares, so will have at most $9$ of these curves covering a square.
Hence since there are only a finite number of curves covering a square,
and since they are images of each other, they must cover a non empty open set
and thus be space filling.
Since these are space filling curves, they
have fractal dimension 2 \cite{Edgar}.

\subsection{Fractal dimension of boundary}
It it known that the fractal
    dimension of the boundary of the fractal dragon curve is
    approximately $1.52$ \cite[p 195]{Edgar}.
    Following the same procedure as Edgar, we now compute the fractal
    dimension of another curve in our family.

    Suppose we apply an infinite sequence of operations, denoted by
    $B=b_0,b_1,\dots$,
    where $b_i\in\{0,1\}$ and $b_i=0$ and $b_i=1$
    indicate application of operations 0
    and $1$ respectively.
    If $B$ is repeating,
    then we will obtain a curve with a self similar nature.
    If there is a finite string before the repeating form,
    ie., $B=ACCCCCC\dots$, then the resulting curve has the form of a finite
    union of curves obtained from $C$, and so the fractal dimension is that
    of $C$.  In the case of non-repeating $B$, then it's not clear how to
    compute the fractal dimension of the boundary.  Now I will describe
    a few of the repeating cases.

        \subsection{Case $0^\infty$ and $1^\infty$}
    As previously mentioned,
    this is the case of Heigway's Fractal dragon curve.
    An example is shown in Figure~\ref{loop}.
    In \cite[p 194,195]{Edgar}
    it is shown that the fractal dimension of the boundary, 
    (Hausdorff and packing and sim dimensions) is approximately $1.52$.

    Since we will use the argument from \cite{Edgar} for several other cases,
    I will reproduce the argument for the computation of the fractal dimension
    of the boundary here, but using our notation.
    Figure~\ref{0000A} shows the first few applications of operation $0$
    to an edge which I label $R$, and which I consider to start from the top.
    I label this line $v$,
    although it probably
    should be labeled $h$, since it crosses a horizontal line, so the images
    in this figure might be better rotated by $90^\circ$.

    Remember that our paths are notated in the form $XvXhXvXhX\dots$
    with $X=L$ or $R$.  Only the $L$ or $R$ are evident in the motion, the
    $v$ and $h$ determine the transformation.
    From Figure~\ref{replace}, the lines we are now using following
    Edgar's construction, have $L$ or $R$ at the end points, and $h$ or $v$
    in the middle, rather than the other way round.
    In this figure, I label the new $L$ and $R$ which appear from the
    $h$ and $v$, in bold to help see how they are introduced.
    The $h$ and $v$ on the segments helps show how although the lines
    look the same at any stage, the $h$ and $v$ mean they will transform
    in different ways.  However, traveling in the opposite direction along
    a line, $R$ and $L$ are reversed, so the left and right sides of the
    path reverse, that is, if we apply an operation $B$ consisting of
    $0$ and $1$s to $L$ and to $R$, then they will have the same boundaries,
    but on opposite sides.  In \cite{Edgar}, the  left side of $0^\infty L$
    is called $U$ and the right side $V$.  I colour these yellow and green
    respectively in the diagram.  I also label points that are in the
    boundary by blue dots.  The yellow and green lines are not in the boundary,
    but they are getting closer and closer to the boundary; the boundaries on either side are the limits of theses lines.  The yellow lines converge to
    scaled copies
    of $U$ and the green converge to scaled copies of $V$.

    The right most image shows the first four
        stages superimposed to help show their relationship to each other.
        The yellow lines indicate the left hand boundary and the green the right hand boundary of $v$ lines, reverse for $h$ lines.  Blue dots corresponding to points lying on the
        boundary of the limit shape.  The two sides of the line should join at the
        dots at either end, but are shown separated for clarity.
        Our aim to to find self similarity, and we are going to express the
        yellow boundary in terms of itself.  This means that once we
        have a yellow component of the left side which is smaller than the
        initial component, we don't have to decompose further, so in
        the third and fourth figure, although we could decompose the
        yellow line from the second step further, we don't bother.
        The yellow horizontal line in the third image is not in the boundary, because
        its end points are not both in the boundary.
        The lower vertical yellow line in this image is in the boundary, but this
        comes from the limiting boundary of the tile to the right of the tile in
        question.
        At this step, the green side, which converges to $V$, is seen to be two
        half size copies of the yellow side $U$, so we don't need to analyse
        further the right side.  In the fourth image, the upper diagonal yellow
        line is a boundary line coming from 
        boundaries
        from the adjacent tiles.  In this image, since we have already decomposed the right side, the colours yellow and green only refer to the left side
        boundary.
        From the fourth picture, we now have
        $U = \frac{\sqrt{2}}{2} U + \frac{\sqrt{2}}{4} U +
        \frac{\sqrt{2}}{4} U $.  This is the self similarity decomposition
        which leads to the computation of the fractal dimension of the boundary
        of the fractal dragon in \cite{Edgar}.
    
    \begin{figure}
      \begin{tikzpicture}[scale={0.1}]
\input{0000A.tex}
      \end{tikzpicture}
      \caption{First steps in construction of $0^\infty L$.
        Operation $0$ applied each time: $h\rightarrow hLv$,
        $v\rightarrow hRv$. }
      \label{0000A}
    \end{figure}

        \subsection{Case $01010101010...$}

\begin{figure}
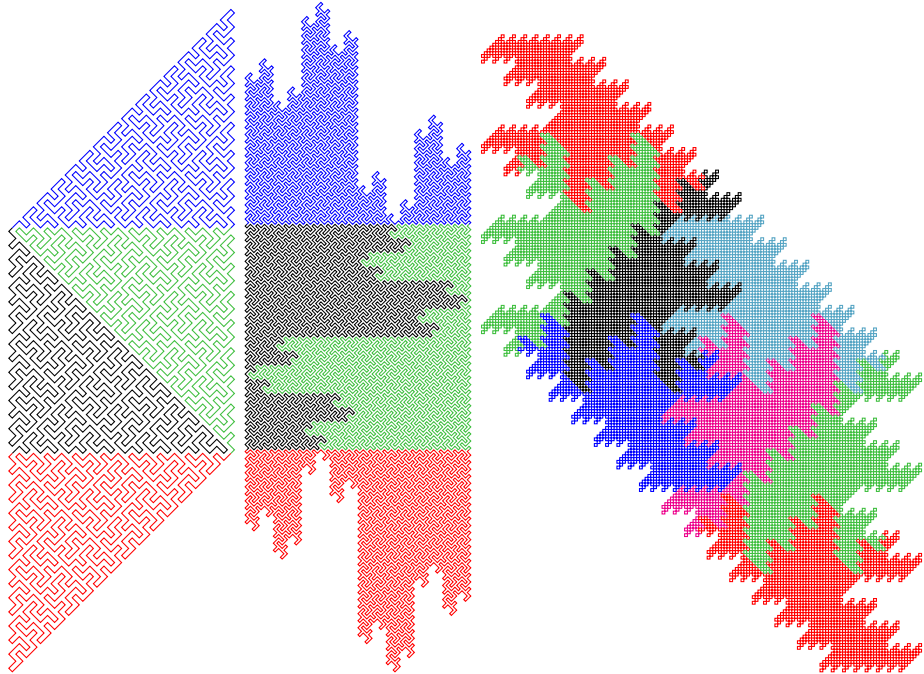

\begin{tikzpicture}
  \input{parallel.tex}
\end{tikzpicture}
\begin{tikzpicture}
  \input{stagger.tex}
  \end{tikzpicture}
\begin{tikzpicture}
\input{thing}
  \end{tikzpicture}
\caption{The closed curves
  $(01)^5LLLL$,
$(0011)^3LLLL$,
  and
  $(001)^3RLLLRLLL$.
}
\label{0101010101}

      \end{figure}

    In this case, a single line ends up converging to a triangle shape,
    and a union of
    four lines in a ring converges to a closed curve with an
    overall parallelogram shape.  The limit of the boundaries is a straight
    line, and so the boundary is not fractal.  Examples of approximations
    are shown in Figure~\ref{0101010101} left.

    \subsection{Case $0011001100110011...$}

    In this case, as in the case of $01010101010...$, a union of two
    of the fractiles forms a square.  However, whilst the boundary is self
    similar and not straight, it appears that the box counting
    dimension is $1$, because the boundary is mostly composed of straight lines.
    This is illustrated in Figure~\ref{0101010101} centre.
    This case appears to be the same
    as one of the examples given in \cite[Figure 7]{Tab}.

    \subsection{Case $(001)^\infty$ and $(110)^\infty$}

    In this case, the fractal dimension of the boundary turns out to be
    approximately $1.4128$.
    We follow \cite[194–195]{Edgar} in making this computation.
    We consider this curve as being produced by repeated application of the
    operation $001$ to the straight line segment corresponding to $Lh$
    The first few applications of this operation are shown in
    Figure~\ref{computefor001recurring}.
    
\comment{    
    \begin{figure}
      \begin{tikzpicture}            
        \input{compute1a.tex}
\begin{scope}[shift={(6.5,0)}]   
  \input{001001001001L.tex}
  \end{scope}
         \end{tikzpicture}
      \caption{Superposition of configurations $L$, $(001)L$ and
        either $(001)^2 L$ (left) or
        $(001)^2 L$ (right). }
    \end{figure}
    }

    \begin{figure}
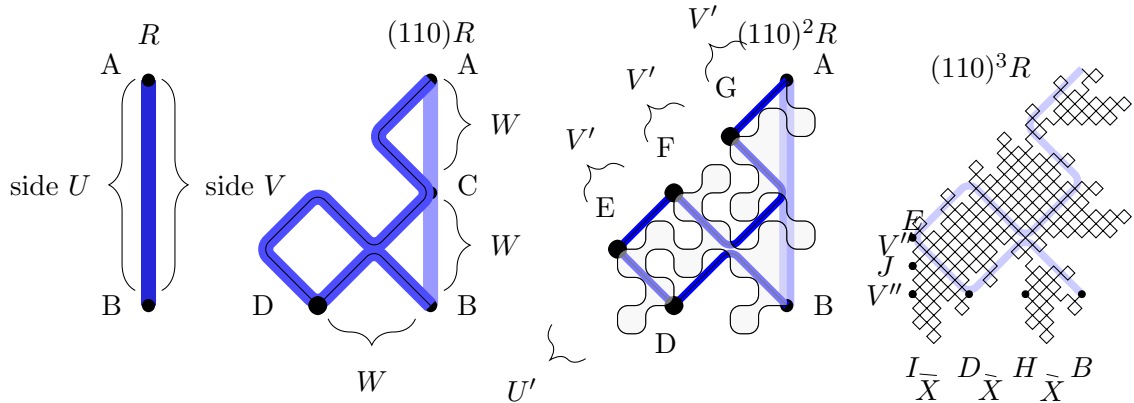

      \begin{tikzpicture}
        \begin{scope}[shift={(-2,0)}]
          \begin{scope}[scale={0.75}]
            \input{compute1b.tex}
          \end{scope}
          \end{scope}          
        \end{tikzpicture}
        \begin{tikzpicture}
          \input{001001001Lb.tex}
      \end{tikzpicture}
      \caption{Self similarity comparison of $(001)^\infty R$.}
\label{computefor001recurring}
    \end{figure}

    Consider Figure~\ref{computefor001recurring}.
    Following Edgar, the left side of $(001)^\infty R$ is labeled $U$
    and the right side $V$.
    Considering $(001)L$, we see that $V$ consists of two identical pieces
    $W$, and that the boundary between points $B$ and $D$ also has form $W$.
    By considering $(001)^2L$, we see that the edge $DE$ is a scaled copy of
    $U$, namely
    $U'= 2^{-3/2}U$, and that the edges $EF, FG, GA$ are all scaled
    copies of the edge $V$, namely $V'=2^{-3/2}V$.
    This results in $U = V' + V' + V' + U' + W$.
    Now considering
    $(001)^4L$ we see (also by comparing with the adjacent tiles, as in
    Figure~\ref{0101010101} right),
    that the edges $W$ have a rotational symmetry about their
    mid points.
    So we now write $W$ as the union of two identical pieces $X$, so
    $W=X + X$.
    Also we now see that the component of the edge from $I$ to $E$ is
    the same as two scaled copies of $V$, namely $V''=2^{-3}V$.
    So we have $U' = X + V'' + V''$, i.e.,
    $2^{-3/2}U = X + 2^{-3}V + 2^{-3}V$, or equivalently,
        $U = 2^{3/2}X + 2^{-3/2}V + 2^{-3/2}V$.
    We also have $V=X + X + X + X$,
    so this gives $U = 2^{3/2}X + 2^{-3/2}(X + X + X + X + X + X)$.
    Now notice that the $X$ edge, as below the $HB$ line in
    is identical, up to a scale factor of $2^{3/2}$, to the edge $BF$ in
    lower left of this line,
    which I will call $BF^{-}$.
From Figure~\ref{computefor001recurring}, we have
$BF^{-}= V' + U' + W$, so together with $BF^{-}=2^{3/2}X$, we have
$2^{3/2}X = V' + U' + W = 2^{-3/2}V + V'' + V'' + X + X + X
= 2^{-3/2}(X + X + X + X) + 2^{-3}(V + V) +  + X + X + X
= 2^{-3/2}(X + X + X + X) + 2^{-3}(X + X+ X + X + X + X+X+X)   + X + X + X
$.
Now we have an expression for $X$ in terms of copies of itself,
namely
$$
X =
2^{-3}(X + X + X + X) + 2^{-9/2}(X + X+ X + X + X + X+X+X)   + 2^{-3/2}(X + X + X).$$
For brevity, I rewrite this formula as
\begin{equation}
X = 2^{-3}X^{\oplus 4} + 2^{-9/2}X^{\oplus 8}   + 2^{-3/2}X^{\oplus 3}
  \label{eqn:sim}
  \end{equation}

Since $U$ and $V$ are both unions of a finite number of scaled copies of $X$,
all the edges have the same fractal dimension as $X$.
To find this dimension, (following \cite[Theorem 6.5.4., p193]{Edgar})
we now need to find $s$, the sim value, such that
$$
4\times (2^{-3})^s + 8 \times (2^{-9/2})^s + 3\times  (2^{-3/2})^s=1.
$$
This implies $2^{-3s/2}$ is a root of
$$8\lambda^3 + 4\lambda^2 + 3\lambda - 1 =0,$$
which has only one real root, $\lambda \approx 0.23017$, so the fractal
dimension is $s=\frac{2}{3}\log_2(1/\lambda)\approx 1.4128$.

\begin{figure}
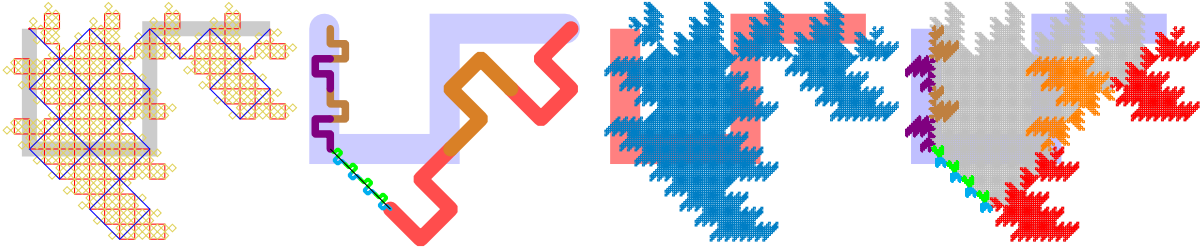

  \begin{tikzpicture}[scale={0.08}]
  \input{001B27septA1.tex}
  \begin{scope}[shift = {(50,0)}]
    \input{Xbreakdown3.tex}
    \end{scope}
    \begin{scope}[shift = {(100,0)}]
      \input{001B27septB1.tex}
    \end{scope}
        \begin{scope}[shift = {(150,0)}]
      \input{001B27septD1.tex}
          \end{scope}
  \end{tikzpicture}
  \caption{A self similar decomposition of $X$ for
    computing the fractal dimension of the boundary of $(001)^\infty L$.
    From left to right, the first image shows superimposed
    $(001)^nLLLR$ for $n=0,\dots, 3$, the second picks out various copies of
    $LLLR$ in the boundary of $(001)^nLLLR$.  The third shows $(001)^4 LLLR$,
    the lower boundary of which is $X$. The fourth shows how applying
    $(001)^\infty$ ($4$ rather than $\infty$ for the purposes of the drawing) to the components in
    the second image results in the lower boundary in the third image,
    showing the self similarity of $X$.
  }
  \label{Xbreakdown}
\end{figure}

Figure~\ref{Xbreakdown} illustrates sets corresponding to
the relationship in Equation~\ref{eqn:sim}.
In this figure,
the blue line, which can be described as $F=LLLR$ is given in blue.
The gray line is $(001)F$, and the black line is $(001)^2F$.
The set $X$ is equal to the boundary of the upper side of
$(001)^\infty F$.
To decompose $X$, we look for
copies of $F$ in the partial iterations $(001)^nX$.

The red and orange lines are the three copies of
$2^{-3/2}F$
The red lines are contained in $(001)F$, and the
orange  is in a neighbour of $(001)F$.
The $4$ copies of  $2^{-3}F$ are given by
the green and cyan lines, which are in
$(001)^2F$ and a neighbouring tile respectively,
and the violet and brown lines are the $8$ copies of
$2^{-9/2}F$.

\subsection{Moran open set condition.}
\begin{figure}
  \input{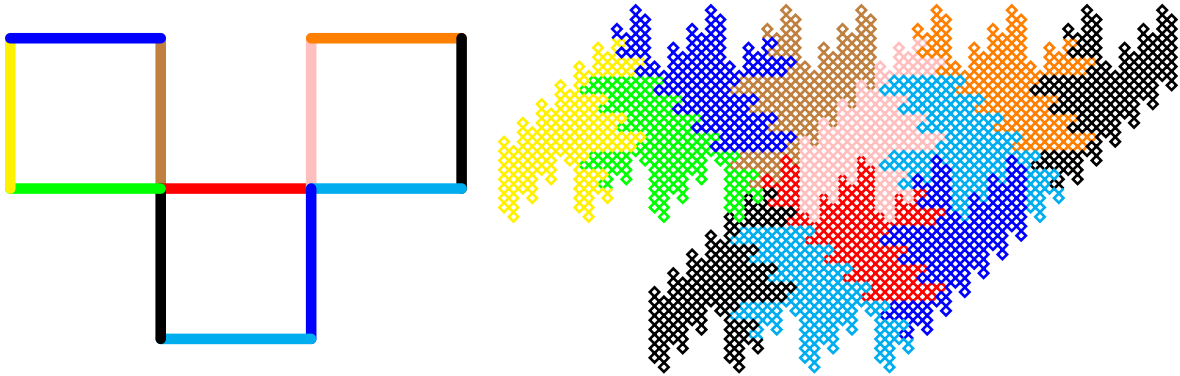}
  \caption{$S=LLLRRRLRRLRR$ and $(001)^3S$.  The interior of the limit
$(001)^\infty S$
    is
    used as a Moran open set.}
  \label{moran}
  \end{figure}

As in \cite{Edgar}, to complete the proof that this is the fractal
dimension, we need to show that the Moran open set condition
\cite[p 191]{Edgar} is satisfied.
The open set condition is that there is an open set $U$ such that for the
iterated system $f_1,f_2,\dots,f_n$, the $f_i(U)$ are disjoint and contained
in $U$.  For us, the $f_i$ are the $15$ functions translating the
blue region in Figure~\ref{Xbreakdown}, third figure,
to the coloured regions round the
lower boundary of this shape in Figure~\ref{Xbreakdown} fourth figure.
A suitable set $U$ is found to be that given in Figure~\ref{moran}.
The left image shows the set $S$ of lines corresponding to the directions
$LLLRRRLRRLRR$ 
The right picture is $(001)^2S$.  Our open set $U$ is the interior of
$(001)^\infty S$, which will look similar to the right figure.

\begin{figure}
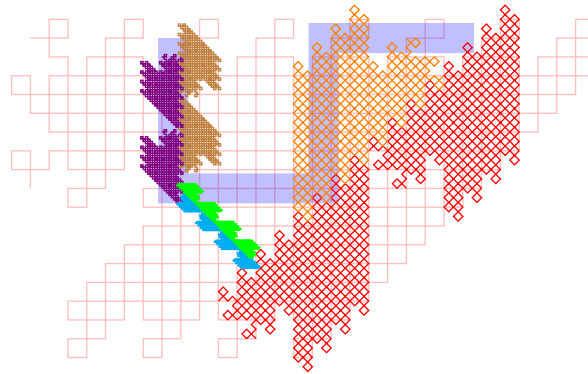

 \begin{tikzpicture}[scale={0.1}]

  \input{manypartsA3.tex}

\begin{scope}[shift={(-100,200)}]
          \bitC{pink}      
\end{scope}
  
    \begin{scope}[shift = {(100,0)}]
           \input{001B27septD3.tex}      
    \end{scope}
 
\end{tikzpicture}
  \caption{Illustration that the Moran open set condition is satisfied.}
  \label{disjointimages}
\end{figure}

I will replace $S$ by a slightly smaller set, where the far left yellow region
is left out, since then it is easier to see that the images $f_i(U)$ are
all disjoint.
Figure~\ref{disjointimages} shows the images of (the slightly smaller)
$U$ under the action of the
$f_i$ illustrating that these regions are all disjoint, and are all contained
in $U$.
In this figure, the pink region is the set $U$, shown in Figure~\ref{moran},
but all in one colour in this figure. On top of the set $U$, all the images
$f_i(U)$ are shown, with the same colours used as in
the fourth image in Figure~\ref{Xbreakdown}.

   \subsection{Case $(0001)^\infty$}
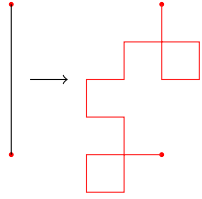
\begin{figure}

  \begin{tikzpicture}[scale={0.1}]
\begin{scope}[shift={(80,-100)}]
\draw[fill=red,red] (0,0) circle (0.25);
\draw[fill=red,red] (0,20) circle (0.25);
\draw(0,0)--(0,20);
\draw[->](2.5,10)--(7.5,10);
\end{scope}
\begin{scope}[rotate={-90}]
\begin{scope}[shift={(100,100)}]
\draw[fill=red,red] (0,0) circle (0.25);
\draw[fill=red,red] (-20,0) circle (0.25);
\draw[red!90!white] (0,0) -- ++(0,-5.000000) -- ++(5,0.000000) -- ++(0,-5.000000) -- ++(-5,0.000000) -- ++(0,5.000000) -- ++(-5,0.000000) -- ++(0,-5.000000) -- ++(-5,0.000000) -- ++(0,5.000000) -- ++(-5,0.000000) -- ++(0,5.000000) -- ++(5,0.000000) -- ++(0,5.000000) -- ++(-5,0.000000) -- ++(0,-5.000000) -- ++(-5,0.000000) ;
\end{scope}
\end{scope}
\end{tikzpicture}
  \caption{The line $L$ and it image $(0001)L$.}
  \label{0001A}
  \end{figure}

\begin{figure}
\begin{tikzpicture}[scale={0.1}]
\begin{scope}[rotate={180}]
\begin{scope}[shift={(200,-200)}]
\filldraw[red] (0,0) circle (0.1);
\draw[blue!50!white, line width = 6] (0,0) -- ++(0,20.000000) 
-- ++(-20,0.000000) -- ++(0,-20.000000) -- ++(-20,0.000000) ;
\end{scope}
\end{scope}

\begin{scope}[rotate={-90}]
\begin{scope}[shift={(-199.99999999999994,-199.99999999999991)}]
\filldraw[red] (0,0) circle (0.1);
\draw[red!50!white, line width=2] (0,0) -- ++(0,5.000000) 
-- ++(-5,0.000000) -- ++(0,5.000000) -- ++(5,0.000000) -- ++(0,-5.000000) -- ++(5,0.000000) -- ++(0,5.000000) -- ++(5,0.000000) -- ++(0,-5.000000) -- ++(5,0.000000) -- ++(0,-5.000000) -- ++(-5,0.000000) -- ++(0,-5.000000) -- ++(5,0.000000) -- ++(0,5.000000) -- ++(5,0.000000) -- ++(0,5.000000) -- ++(-5,0.000000) -- ++(0,5.000000) -- ++(5,0.000000) -- ++(0,-5.000000) -- ++(5,0.000000) -- ++(0,5.000000) -- ++(5,0.000000) -- ++(0,5.000000) -- ++(-5,0.000000) -- ++(0,5.000000) -- ++(5,0.000000) -- ++(0,5.000000) -- ++(-5,0.000000) -- ++(0,-5.000000) -- ++(-5,0.000000) -- ++(0,-5.000000) -- ++(5,0.000000) -- ++(0,-5.000000) -- ++(-5,0.000000) -- ++(0,5.000000) -- ++(-5,0.000000) -- ++(0,-5.000000) -- ++(-5,0.000000) -- ++(0,5.000000) -- ++(-5,0.000000) -- ++(0,5.000000) -- ++(5,0.000000) -- ++(0,5.000000) -- ++(-5,0.000000) -- ++(0,-5.000000) -- ++(-5,0.000000) -- ++(0,5.000000) -- ++(-5,0.000000) -- ++(0,5.000000) -- ++(5,0.000000) -- ++(0,-5.000000) -- ++(5,0.000000) -- ++(0,5.000000) -- ++(5,0.000000) -- ++(0,5.000000) -- ++(-5,0.000000) -- ++(0,5.000000) -- ++(5,0.000000) -- ++(0,5.000000) -- ++(-5,0.000000) -- ++(0,-5.000000) -- ++(-5,0.000000) ;
\end{scope}
\end{scope}

\begin{scope}[rotate={0}]
\begin{scope}[shift={(-199.99999999999994,199.99999999999994)}]
\draw[black] (0,0) -- ++(0,1.250000) 
-- ++(-1.25,0.000000) -- ++(0,1.250000) -- ++(1.25,0.000000) -- ++(0,-1.250000) -- ++(1.25,0.000000) -- ++(0,1.250000) -- ++(1.25,0.000000) -- ++(0,-1.250000) -- ++(1.25,0.000000) -- ++(0,-1.250000) -- ++(-1.25,0.000000) -- ++(0,-1.250000) -- ++(1.25,0.000000) -- ++(0,1.250000) -- ++(1.25,0.000000) -- ++(0,1.250000) -- ++(-1.25,0.000000) -- ++(0,1.250000) -- ++(1.25,0.000000) -- ++(0,-1.250000) -- ++(1.25,0.000000) -- ++(0,1.250000) -- ++(1.25,0.000000) -- ++(0,1.250000) -- ++(-1.25,0.000000) -- ++(0,1.250000) -- ++(1.25,0.000000) -- ++(0,1.250000) -- ++(-1.25,0.000000) -- ++(0,-1.250000) -- ++(-1.25,0.000000) -- ++(0,1.250000) -- ++(-1.25,0.000000) -- ++(0,1.250000) -- ++(1.25,0.000000) -- ++(0,-1.250000) -- ++(1.25,0.000000) -- ++(0,1.250000) -- ++(1.25,0.000000) -- ++(0,-1.250000) -- ++(1.25,0.000000) -- ++(0,-1.250000) -- ++(-1.25,0.000000) -- ++(0,-1.250000) -- ++(1.25,0.000000) -- ++(0,1.250000) -- ++(1.25,0.000000) -- ++(0,-1.250000) -- ++(1.25,0.000000) -- ++(0,-1.250000) -- ++(-1.25,0.000000) -- ++(0,1.250000) -- ++(-1.25,0.000000) -- ++(0,-1.250000) -- ++(-1.25,0.000000) -- ++(0,-1.250000) -- ++(1.25,0.000000) -- ++(0,-1.250000) -- ++(-1.25,0.000000) -- ++(0,-1.250000) -- ++(1.25,0.000000) -- ++(0,1.250000) -- ++(1.25,0.000000) -- ++(0,-1.250000) -- ++(1.25,0.000000) -- ++(0,-1.250000) -- ++(-1.25,0.000000) -- ++(0,1.250000) -- ++(-1.25,0.000000) -- ++(0,-1.250000) -- ++(-1.25,0.000000) -- ++(0,1.250000) -- ++(-1.25,0.000000) -- ++(0,1.250000) -- ++(1.25,0.000000) -- ++(0,1.250000) -- ++(-1.25,0.000000) -- ++(0,-1.250000) -- ++(-1.25,0.000000) -- ++(0,-1.250000) -- ++(1.25,0.000000) -- ++(0,-1.250000) -- ++(-1.25,0.000000) -- ++(0,1.250000) -- ++(-1.25,0.000000) -- ++(0,-1.250000) -- ++(-1.25,0.000000) -- ++(0,-1.250000) -- ++(1.25,0.000000) -- ++(0,-1.250000) -- ++(-1.25,0.000000) -- ++(0,-1.250000) -- ++(1.25,0.000000) -- ++(0,1.250000) -- ++(1.25,0.000000) -- ++(0,1.250000) -- ++(-1.25,0.000000) -- ++(0,1.250000) -- ++(1.25,0.000000) -- ++(0,-1.250000) -- ++(1.25,0.000000) -- ++(0,1.250000) -- ++(1.25,0.000000) -- ++(0,-1.250000) -- ++(1.25,0.000000) -- ++(0,-1.250000) -- ++(-1.25,0.000000) -- ++(0,-1.250000) -- ++(1.25,0.000000) -- ++(0,1.250000) -- ++(1.25,0.000000) -- ++(0,-1.250000) -- ++(1.25,0.000000) -- ++(0,-1.250000) -- ++(-1.25,0.000000) -- ++(0,1.250000) -- ++(-1.25,0.000000) -- ++(0,-1.250000) -- ++(-1.25,0.000000) -- ++(0,-1.250000) -- ++(1.25,0.000000) -- ++(0,-1.250000) -- ++(-1.25,0.000000) -- ++(0,-1.250000) -- ++(1.25,0.000000) -- ++(0,1.250000) -- ++(1.25,0.000000) -- ++(0,-1.250000) -- ++(1.25,0.000000) -- ++(0,-1.250000) -- ++(-1.25,0.000000) -- ++(0,1.250000) -- ++(-1.25,0.000000) -- ++(0,-1.250000) -- ++(-1.25,0.000000) -- ++(0,1.250000) -- ++(-1.25,0.000000) -- ++(0,1.250000) -- ++(1.25,0.000000) -- ++(0,1.250000) -- ++(-1.25,0.000000) -- ++(0,-1.250000) -- ++(-1.25,0.000000) -- ++(0,-1.250000) -- ++(1.25,0.000000) -- ++(0,-1.250000) -- ++(-1.25,0.000000) -- ++(0,1.250000) -- ++(-1.25,0.000000) -- ++(0,-1.250000) 
-- ++(-1.25,0.000000) -- ++(0,-1.250000) -- ++(1.25,0.000000) -- ++(0,-1.250000) -- ++(-1.25,0.000000) -- ++(0,-1.250000) -- ++(1.25,0.000000) -- ++(0,1.250000) -- ++(1.25,0.000000) -- ++(0,-1.250000) -- ++(1.25,0.000000) -- ++(0,-1.250000) -- ++(-1.25,0.000000) -- ++(0,1.250000) -- ++(-1.25,0.000000) -- ++(0,-1.250000) -- ++(-1.25,0.000000) -- ++(0,1.250000) -- ++(-1.25,0.000000) -- ++(0,1.250000) -- ++(1.25,0.000000) -- ++(0,1.250000) -- ++(-1.25,0.000000) -- ++(0,-1.250000) -- ++(-1.25,0.000000) -- ++(0,1.250000) -- ++(-1.25,0.000000) -- ++(0,1.250000) -- ++(1.25,0.000000) -- ++(0,-1.250000) -- ++(1.25,0.000000) -- ++(0,1.250000) -- ++(1.25,0.000000) -- ++(0,1.250000) -- ++(-1.25,0.000000) -- ++(0,1.250000) -- ++(1.25,0.000000) -- ++(0,1.250000) -- ++(-1.25,0.000000) -- ++(0,-1.250000) -- ++(-1.25,0.000000) -- ++(0,-1.250000) -- ++(1.25,0.000000) -- ++(0,-1.250000) -- ++(-1.25,0.000000) -- ++(0,1.250000) -- ++(-1.25,0.000000) -- ++(0,-1.250000) -- ++(-1.25,0.000000) -- ++(0,1.250000) -- ++(-1.25,0.000000) -- ++(0,1.250000) -- ++(1.25,0.000000) -- ++(0,1.250000) -- ++(-1.25,0.000000) -- ++(0,-1.250000) -- ++(-1.25,0.000000) -- ++(0,-1.250000) -- ++(1.25,0.000000) -- ++(0,-1.250000) -- ++(-1.25,0.000000) -- ++(0,1.250000) -- ++(-1.25,0.000000) -- ++(0,-1.250000) -- ++(-1.25,0.000000) -- ++(0,-1.250000) -- ++(1.25,0.000000) -- ++(0,-1.250000) -- ++(-1.25,0.000000) -- ++(0,-1.250000) -- ++(1.25,0.000000) -- ++(0,1.250000) -- ++(1.25,0.000000) -- ++(0,1.250000) -- ++(-1.25,0.000000) -- ++(0,1.250000) -- ++(1.25,0.000000) -- ++(0,-1.250000) -- ++(1.25,0.000000) -- ++(0,1.250000) -- ++(1.25,0.000000) -- ++(0,-1.250000) -- ++(1.25,0.000000) -- ++(0,-1.250000) -- ++(-1.25,0.000000) -- ++(0,-1.250000) -- ++(1.25,0.000000) -- ++(0,1.250000) -- ++(1.25,0.000000) -- ++(0,-1.250000) -- ++(1.25,0.000000) -- ++(0,-1.250000) -- ++(-1.25,0.000000) -- ++(0,1.250000) -- ++(-1.25,0.000000) -- ++(0,-1.250000) -- ++(-1.25,0.000000) -- ++(0,-1.250000) -- ++(1.25,0.000000) -- ++(0,-1.250000) -- ++(-1.25,0.000000) -- ++(0,-1.250000) -- ++(1.25,0.000000) -- ++(0,1.250000) -- ++(1.25,0.000000) -- ++(0,1.250000) -- ++(-1.25,0.000000) -- ++(0,1.250000) -- ++(1.25,0.000000) -- ++(0,-1.250000) -- ++(1.25,0.000000) -- ++(0,1.250000) -- ++(1.25,0.000000) -- ++(0,-1.250000) -- ++(1.25,0.000000) -- ++(0,-1.250000) -- ++(-1.25,0.000000) -- ++(0,-1.250000) -- ++(1.25,0.000000) -- ++(0,1.250000) -- ++(1.25,0.000000) -- ++(0,1.250000) -- ++(-1.25,0.000000) -- ++(0,1.250000) -- ++(1.25,0.000000) -- ++(0,-1.250000) -- ++(1.25,0.000000) -- ++(0,1.250000) -- ++(1.25,0.000000) -- ++(0,1.250000) -- ++(-1.25,0.000000) -- ++(0,1.250000) -- ++(1.25,0.000000) -- ++(0,1.250000) -- ++(-1.25,0.000000) -- ++(0,-1.250000) -- ++(-1.25,0.000000) -- ++(0,1.250000) -- ++(-1.25,0.000000) -- ++(0,1.250000) -- ++(1.25,0.000000) -- ++(0,-1.250000) -- ++(1.25,0.000000) -- ++(0,1.250000) -- ++(1.25,0.000000) -- ++(0,-1.250000) -- ++(1.25,0.000000) -- ++(0,-1.250000) -- ++(-1.25,0.000000) -- ++(0,-1.250000) 
-- ++(1.25,0.000000) -- ++(0,1.250000) -- ++(1.25,0.000000) -- ++(0,-1.250000) -- ++(1.25,0.000000) -- ++(0,-1.250000) -- ++(-1.25,0.000000) -- ++(0,1.250000) -- ++(-1.25,0.000000) -- ++(0,-1.250000) -- ++(-1.25,0.000000) -- ++(0,-1.250000) -- ++(1.25,0.000000) -- ++(0,-1.250000) -- ++(-1.25,0.000000) -- ++(0,-1.250000) -- ++(1.25,0.000000) -- ++(0,1.250000) -- ++(1.25,0.000000) -- ++(0,-1.250000) -- ++(1.25,0.000000) -- ++(0,-1.250000) -- ++(-1.25,0.000000) -- ++(0,1.250000) -- ++(-1.25,0.000000) -- ++(0,-1.250000) -- ++(-1.25,0.000000) -- ++(0,1.250000) -- ++(-1.25,0.000000) -- ++(0,1.250000) -- ++(1.25,0.000000) -- ++(0,1.250000) -- ++(-1.25,0.000000) -- ++(0,-1.250000) -- ++(-1.25,0.000000) -- ++(0,-1.250000) -- ++(1.25,0.000000) -- ++(0,-1.250000) -- ++(-1.25,0.000000) -- ++(0,1.250000) -- ++(-1.25,0.000000) -- ++(0,-1.250000) -- ++(-1.25,0.000000) -- ++(0,-1.250000) -- ++(1.25,0.000000) -- ++(0,-1.250000) -- ++(-1.25,0.000000) -- ++(0,-1.250000) -- ++(1.25,0.000000) -- ++(0,1.250000) 
-- ++(1.25,0.000000) -- ++(0,1.250000) -- ++(-1.25,0.000000) -- ++(0,1.250000) -- ++(1.25,0.000000) -- ++(0,-1.250000) -- ++(1.25,0.000000) -- ++(0,1.250000) -- ++(1.25,0.000000) -- ++(0,-1.250000) -- ++(1.25,0.000000) -- ++(0,-1.250000) -- ++(-1.25,0.000000) -- ++(0,-1.250000) -- ++(1.25,0.000000) -- ++(0,1.250000) -- ++(1.25,0.000000) -- ++(0,-1.250000) -- ++(1.25,0.000000) -- ++(0,-1.250000) -- ++(-1.25,0.000000) -- ++(0,1.250000) -- ++(-1.25,0.000000) -- ++(0,-1.250000) -- ++(-1.25,0.000000) -- ++(0,-1.250000) -- ++(1.25,0.000000) -- ++(0,-1.250000) -- ++(-1.25,0.000000) -- ++(0,-1.250000) -- ++(1.25,0.000000) -- ++(0,1.250000) -- ++(1.25,0.000000) -- ++(0,1.250000) -- ++(-1.25,0.000000) -- ++(0,1.250000) -- ++(1.25,0.000000) -- ++(0,-1.250000) -- ++(1.25,0.000000) -- ++(0,1.250000) -- ++(1.25,0.000000) -- ++(0,-1.250000) -- ++(1.25,0.000000) -- ++(0,-1.250000) -- ++(-1.25,0.000000) -- ++(0,-1.250000) -- ++(1.25,0.000000) -- ++(0,1.250000) -- ++(1.25,0.000000) -- ++(0,1.250000) 
-- ++(-1.25,0.000000) -- ++(0,1.250000) -- ++(1.25,0.000000) -- ++(0,-1.250000) -- ++(1.25,0.000000) -- ++(0,1.250000) -- ++(1.25,0.000000) -- ++(0,1.250000) -- ++(-1.25,0.000000) -- ++(0,1.250000) -- ++(1.25,0.000000) -- ++(0,1.250000) -- ++(-1.25,0.000000) -- ++(0,-1.250000) -- ++(-1.25,0.000000) -- ++(0,1.250000) -- ++(-1.25,0.000000) -- ++(0,1.250000) -- ++(1.25,0.000000) -- ++(0,-1.250000) -- ++(1.25,0.000000) -- ++(0,1.250000) -- ++(1.25,0.000000) -- ++(0,-1.250000) -- ++(1.25,0.000000) -- ++(0,-1.250000) -- ++(-1.25,0.000000) -- ++(0,-1.250000) -- ++(1.25,0.000000) -- ++(0,1.250000) -- ++(1.25,0.000000) -- ++(0,-1.250000) -- ++(1.25,0.000000) -- ++(0,-1.250000) -- ++(-1.25,0.000000) -- ++(0,1.250000) -- ++(-1.25,0.000000) -- ++(0,-1.250000) -- ++(-1.25,0.000000) -- ++(0,-1.250000) -- ++(1.25,0.000000) -- ++(0,-1.250000) -- ++(-1.25,0.000000) -- ++(0,-1.250000) -- ++(1.25,0.000000) -- ++(0,1.250000) -- ++(1.25,0.000000) -- ++(0,1.250000) -- ++(-1.25,0.000000) -- ++(0,1.250000) -- ++(1.25,0.000000) -- ++(0,-1.250000) -- ++(1.25,0.000000) -- ++(0,1.250000) -- ++(1.25,0.000000) -- ++(0,-1.250000) -- ++(1.25,0.000000) -- ++(0,-1.250000) -- ++(-1.25,0.000000) -- ++(0,-1.250000) -- ++(1.25,0.000000) -- ++(0,1.250000) -- ++(1.25,0.000000) -- ++(0,1.250000) -- ++(-1.25,0.000000) -- ++(0,1.250000) -- ++(1.25,0.000000) -- ++(0,-1.250000) -- ++(1.25,0.000000) -- ++(0,1.250000) -- ++(1.25,0.000000) -- ++(0,1.250000) -- ++(-1.25,0.000000) -- ++(0,1.250000) -- ++(1.25,0.000000) -- ++(0,1.250000) -- ++(-1.25,0.000000) -- ++(0,-1.250000) -- ++(-1.25,0.000000) -- ++(0,-1.250000) -- ++(1.25,0.000000) -- ++(0,-1.250000) -- ++(-1.25,0.000000) -- ++(0,1.250000) -- ++(-1.25,0.000000) -- ++(0,-1.250000) -- ++(-1.25,0.000000) -- ++(0,1.250000) -- ++(-1.25,0.000000) -- ++(0,1.250000) -- ++(1.25,0.000000) -- ++(0,1.250000) -- ++(-1.25,0.000000) -- ++(0,-1.250000) -- ++(-1.25,0.000000) -- ++(0,1.250000) -- ++(-1.25,0.000000) -- ++(0,1.250000) -- ++(1.25,0.000000) -- ++(0,-1.250000) -- ++(1.25,0.000000) -- ++(0,1.250000) -- ++(1.25,0.000000) -- ++(0,1.250000) -- ++(-1.25,0.000000) -- ++(0,1.250000) -- ++(1.25,0.000000) -- ++(0,1.250000) -- ++(-1.25,0.000000) -- ++(0,-1.250000) -- ++(-1.25,0.000000) -- ++(0,-1.250000) -- ++(1.25,0.000000) -- ++(0,-1.250000) -- ++(-1.25,0.000000) -- ++(0,1.250000) -- ++(-1.25,0.000000) -- ++(0,-1.250000) -- ++(-1.25,0.000000) -- ++(0,1.250000) -- ++(-1.25,0.000000) -- ++(0,1.250000) -- ++(1.25,0.000000) -- ++(0,1.250000) -- ++(-1.25,0.000000) -- ++(0,-1.250000) -- ++(-1.25,0.000000) -- ++(0,-1.250000) -- ++(1.25,0.000000) -- ++(0,-1.250000) -- ++(-1.25,0.000000) -- ++(0,1.250000) -- ++(-1.25,0.000000) -- ++(0,-1.250000) -- ++(-1.25,0.000000) -- ++(0,-1.250000) -- ++(1.25,0.000000) -- ++(0,-1.250000) -- ++(-1.25,0.000000) -- ++(0,-1.250000) -- ++(1.25,0.000000) -- ++(0,1.250000) -- ++(1.25,0.000000) -- ++(0,-1.250000) -- ++(1.25,0.000000) -- ++(0,-1.250000) -- ++(-1.25,0.000000) -- ++(0,1.250000) -- ++(-1.25,0.000000) -- ++(0,-1.250000) 
-- ++(-1.25,0.000000) -- ++(0,1.250000) -- ++(-1.25,0.000000) -- ++(0,1.250000) -- ++(1.25,0.000000) -- ++(0,1.250000) -- ++(-1.25,0.000000) -- ++(0,-1.250000) -- ++(-1.25,0.000000) -- ++(0,1.250000) -- ++(-1.25,0.000000) -- ++(0,1.250000) -- ++(1.25,0.000000) -- ++(0,-1.250000) -- ++(1.25,0.000000) -- ++(0,1.250000) -- ++(1.25,0.000000) -- ++(0,1.250000) -- ++(-1.25,0.000000) -- ++(0,1.250000) -- ++(1.25,0.000000) -- ++(0,1.250000) -- ++(-1.25,0.000000) -- ++(0,-1.250000) -- ++(-1.25,0.000000) -- ++(0,1.250000) -- ++(-1.25,0.000000) -- ++(0,1.250000) -- ++(1.25,0.000000) -- ++(0,-1.250000) -- ++(1.25,0.000000) -- ++(0,1.250000) -- ++(1.25,0.000000) -- ++(0,-1.250000) -- ++(1.25,0.000000) -- ++(0,-1.250000) -- ++(-1.25,0.000000) -- ++(0,-1.250000) -- ++(1.25,0.000000) -- ++(0,1.250000) -- ++(1.25,0.000000) -- ++(0,1.250000) -- ++(-1.25,0.000000) -- ++(0,1.250000) -- ++(1.25,0.000000) -- ++(0,-1.250000) -- ++(1.25,0.000000) -- ++(0,1.250000) -- ++(1.25,0.000000) -- ++(0,1.250000) 
-- ++(-1.25,0.000000) -- ++(0,1.250000) -- ++(1.25,0.000000) -- ++(0,1.250000) -- ++(-1.25,0.000000) -- ++(0,-1.250000) -- ++(-1.25,0.000000) -- ++(0,-1.250000) -- ++(1.25,0.000000) -- ++(0,-1.250000) -- ++(-1.25,0.000000) -- ++(0,1.250000) -- ++(-1.25,0.000000) -- ++(0,-1.250000) -- ++(-1.25,0.000000) -- ++(0,1.250000) -- ++(-1.25,0.000000) -- ++(0,1.250000) -- ++(1.25,0.000000) -- ++(0,1.250000) -- ++(-1.25,0.000000) -- ++(0,-1.250000) -- ++(-1.25,0.000000) -- ++(0,1.250000) -- ++(-1.25,0.000000) -- ++(0,1.250000) -- ++(1.25,0.000000) -- ++(0,-1.250000) -- ++(1.25,0.000000) -- ++(0,1.250000) -- ++(1.25,0.000000) -- ++(0,1.250000) -- ++(-1.25,0.000000) -- ++(0,1.250000) -- ++(1.25,0.000000) -- ++(0,1.250000) -- ++(-1.25,0.000000) -- ++(0,-1.250000) -- ++(-1.25,0.000000) -- ++(0,1.250000) -- ++(-1.25,0.000000) -- ++(0,1.250000) -- ++(1.25,0.000000) -- ++(0,-1.250000) -- ++(1.25,0.000000) -- ++(0,1.250000) -- ++(1.25,0.000000) -- ++(0,-1.250000) -- ++(1.25,0.000000) -- ++(0,-1.250000) -- ++(-1.25,0.000000) -- ++(0,-1.250000) -- ++(1.25,0.000000) -- ++(0,1.250000) -- ++(1.25,0.000000) -- ++(0,1.250000) -- ++(-1.25,0.000000) -- ++(0,1.250000) -- ++(1.25,0.000000) -- ++(0,-1.250000) -- ++(1.25,0.000000) -- ++(0,1.250000) -- ++(1.25,0.000000) -- ++(0,1.250000) -- ++(-1.25,0.000000) -- ++(0,1.250000) -- ++(1.25,0.000000) -- ++(0,1.250000) -- ++(-1.25,0.000000) -- ++(0,-1.250000) -- ++(-1.25,0.000000) -- ++(0,1.250000) -- ++(-1.25,0.000000) -- ++(0,1.250000) -- ++(1.25,0.000000) -- ++(0,-1.250000) -- ++(1.25,0.000000) -- ++(0,1.250000) -- ++(1.25,0.000000) -- ++(0,-1.250000) -- ++(1.25,0.000000) -- ++(0,-1.250000) -- ++(-1.25,0.000000) -- ++(0,-1.250000) -- ++(1.25,0.000000) -- ++(0,1.250000) -- ++(1.25,0.000000) -- ++(0,-1.250000) -- ++(1.25,0.000000) -- ++(0,-1.250000) -- ++(-1.25,0.000000) -- ++(0,1.250000) -- ++(-1.25,0.000000) -- ++(0,-1.250000) -- ++(-1.25,0.000000) -- ++(0,-1.250000) -- ++(1.25,0.000000) -- ++(0,-1.250000) -- ++(-1.25,0.000000) -- ++(0,-1.250000) 
-- ++(1.25,0.000000) -- ++(0,1.250000) -- ++(1.25,0.000000) -- ++(0,1.250000) -- ++(-1.25,0.000000) -- ++(0,1.250000) -- ++(1.25,0.000000) -- ++(0,-1.250000) -- ++(1.25,0.000000) -- ++(0,1.250000) -- ++(1.25,0.000000) -- ++(0,-1.250000) -- ++(1.25,0.000000) -- ++(0,-1.250000) -- ++(-1.25,0.000000) -- ++(0,-1.250000) -- ++(1.25,0.000000) -- ++(0,1.250000) -- ++(1.25,0.000000) -- ++(0,1.250000) -- ++(-1.25,0.000000) -- ++(0,1.250000) -- ++(1.25,0.000000) -- ++(0,-1.250000) -- ++(1.25,0.000000) -- ++(0,1.250000) -- ++(1.25,0.000000) -- ++(0,1.250000) -- ++(-1.25,0.000000) -- ++(0,1.250000) -- ++(1.25,0.000000) -- ++(0,1.250000) -- ++(-1.25,0.000000) -- ++(0,-1.250000) -- ++(-1.25,0.000000) -- ++(0,-1.250000) -- ++(1.25,0.000000) -- ++(0,-1.250000) -- ++(-1.25,0.000000) -- ++(0,1.250000) -- ++(-1.25,0.000000) -- ++(0,-1.250000) -- ++(-1.25,0.000000) -- ++(0,1.250000) -- ++(-1.25,0.000000) -- ++(0,1.250000) -- ++(1.25,0.000000) -- ++(0,1.250000) -- ++(-1.25,0.000000) -- ++(0,-1.250000) 
-- ++(-1.25,0.000000) -- ++(0,1.250000) -- ++(-1.25,0.000000) -- ++(0,1.250000) -- ++(1.25,0.000000) -- ++(0,-1.250000) -- ++(1.25,0.000000) -- ++(0,1.250000) -- ++(1.25,0.000000) -- ++(0,1.250000) -- ++(-1.25,0.000000) -- ++(0,1.250000) -- ++(1.25,0.000000) -- ++(0,1.250000) -- ++(-1.25,0.000000) -- ++(0,-1.250000) -- ++(-1.25,0.000000) -- ++(0,1.250000) -- ++(-1.25,0.000000) -- ++(0,1.250000) -- ++(1.25,0.000000) -- ++(0,-1.250000) -- ++(1.25,0.000000) -- ++(0,1.250000) -- ++(1.25,0.000000) -- ++(0,-1.250000) -- ++(1.25,0.000000) -- ++(0,-1.250000) -- ++(-1.25,0.000000) -- ++(0,-1.250000) -- ++(1.25,0.000000) -- ++(0,1.250000) -- ++(1.25,0.000000) -- ++(0,1.250000) -- ++(-1.25,0.000000) -- ++(0,1.250000) -- ++(1.25,0.000000) -- ++(0,-1.250000) -- ++(1.25,0.000000) -- ++(0,1.250000) -- ++(1.25,0.000000) -- ++(0,1.250000) -- ++(-1.25,0.000000) -- ++(0,1.250000) -- ++(1.25,0.000000) -- ++(0,1.250000) -- ++(-1.25,0.000000) -- ++(0,-1.250000) -- ++(-1.25,0.000000) -- ++(0,1.250000) -- ++(-1.25,0.000000) -- ++(0,1.250000) -- ++(1.25,0.000000) -- ++(0,-1.250000) -- ++(1.25,0.000000) -- ++(0,1.250000) -- ++(1.25,0.000000) -- ++(0,-1.250000) -- ++(1.25,0.000000) -- ++(0,-1.250000) -- ++(-1.25,0.000000) -- ++(0,-1.250000) -- ++(1.25,0.000000) -- ++(0,1.250000) -- ++(1.25,0.000000) -- ++(0,-1.250000) -- ++(1.25,0.000000) -- ++(0,-1.250000) -- ++(-1.25,0.000000) -- ++(0,1.250000) -- ++(-1.25,0.000000) -- ++(0,-1.250000) -- ++(-1.25,0.000000) -- ++(0,-1.250000) -- ++(1.25,0.000000) -- ++(0,-1.250000) -- ++(-1.25,0.000000) -- ++(0,-1.250000) -- ++(1.25,0.000000) -- ++(0,1.250000) -- ++(1.25,0.000000) -- ++(0,-1.250000) -- ++(1.25,0.000000) -- ++(0,-1.250000) -- ++(-1.25,0.000000) -- ++(0,1.250000) -- ++(-1.25,0.000000) -- ++(0,-1.250000) -- ++(-1.25,0.000000) -- ++(0,1.250000) -- ++(-1.25,0.000000) -- ++(0,1.250000) -- ++(1.25,0.000000) -- ++(0,1.250000) -- ++(-1.25,0.000000) -- ++(0,-1.250000) -- ++(-1.25,0.000000) -- ++(0,-1.250000) -- ++(1.25,0.000000) -- ++(0,-1.250000) -- ++(-1.25,0.000000) -- ++(0,1.250000) -- ++(-1.25,0.000000) -- ++(0,-1.250000) -- ++(-1.25,0.000000) -- ++(0,-1.250000) -- ++(1.25,0.000000) -- ++(0,-1.250000) -- ++(-1.25,0.000000) -- ++(0,-1.250000) -- ++(1.25,0.000000) -- ++(0,1.250000) -- ++(1.25,0.000000) -- ++(0,1.250000) -- ++(-1.25,0.000000) -- ++(0,1.250000) -- ++(1.25,0.000000) -- ++(0,-1.250000) -- ++(1.25,0.000000) -- ++(0,1.250000) -- ++(1.25,0.000000) -- ++(0,-1.250000) -- ++(1.25,0.000000) -- ++(0,-1.250000) -- ++(-1.25,0.000000) -- ++(0,-1.250000) -- ++(1.25,0.000000) -- ++(0,1.250000) -- ++(1.25,0.000000) -- ++(0,-1.250000) -- ++(1.25,0.000000) -- ++(0,-1.250000) -- ++(-1.25,0.000000) -- ++(0,1.250000) -- ++(-1.25,0.000000) -- ++(0,-1.250000) -- ++(-1.25,0.000000) -- ++(0,-1.250000) -- ++(1.25,0.000000) -- ++(0,-1.250000) -- ++(-1.25,0.000000) -- ++(0,-1.250000) -- ++(1.25,0.000000) -- ++(0,1.250000) -- ++(1.25,0.000000) -- ++(0,1.250000) -- ++(-1.25,0.000000) -- ++(0,1.250000) -- ++(1.25,0.000000) -- ++(0,-1.250000) -- ++(1.25,0.000000) -- ++(0,1.250000) -- ++(1.25,0.000000) -- ++(0,-1.250000) -- ++(1.25,0.000000) -- ++(0,-1.250000) -- ++(-1.25,0.000000) -- ++(0,-1.250000) -- ++(1.25,0.000000) -- ++(0,1.250000) -- ++(1.25,0.000000) -- ++(0,1.250000) -- ++(-1.25,0.000000) -- ++(0,1.250000) -- ++(1.25,0.000000) -- ++(0,-1.250000) -- ++(1.25,0.000000) -- ++(0,1.250000) -- ++(1.25,0.000000) -- ++(0,1.250000) -- ++(-1.25,0.000000) -- ++(0,1.250000) -- ++(1.25,0.000000) -- ++(0,1.250000) -- ++(-1.25,0.000000) -- ++(0,-1.250000) -- ++(-1.25,0.000000) -- ++(0,1.250000) -- ++(-1.25,0.000000) -- ++(0,1.250000) -- ++(1.25,0.000000) -- ++(0,-1.250000) -- ++(1.25,0.000000) -- ++(0,1.250000) -- ++(1.25,0.000000) -- ++(0,-1.250000) -- ++(1.25,0.000000) -- ++(0,-1.250000) -- ++(-1.25,0.000000) -- ++(0,-1.250000) -- ++(1.25,0.000000) -- ++(0,1.250000) -- ++(1.25,0.000000) -- ++(0,-1.250000) -- ++(1.25,0.000000) -- ++(0,-1.250000) -- ++(-1.25,0.000000) -- ++(0,1.250000) -- ++(-1.25,0.000000) -- ++(0,-1.250000) 
-- ++(-1.25,0.000000) -- ++(0,-1.250000) -- ++(1.25,0.000000) -- ++(0,-1.250000) -- ++(-1.25,0.000000) -- ++(0,-1.250000) -- ++(1.25,0.000000) -- ++(0,1.250000) -- ++(1.25,0.000000) -- ++(0,1.250000) -- ++(-1.25,0.000000) -- ++(0,1.250000) -- ++(1.25,0.000000) -- ++(0,-1.250000) -- ++(1.25,0.000000) -- ++(0,1.250000) -- ++(1.25,0.000000) -- ++(0,-1.250000) -- ++(1.25,0.000000) -- ++(0,-1.250000) -- ++(-1.25,0.000000) -- ++(0,-1.250000) -- ++(1.25,0.000000) -- ++(0,1.250000) -- ++(1.25,0.000000) -- ++(0,1.250000) -- ++(-1.25,0.000000) -- ++(0,1.250000) -- ++(1.25,0.000000) -- ++(0,-1.250000) -- ++(1.25,0.000000) -- ++(0,1.250000) -- ++(1.25,0.000000) -- ++(0,1.250000) -- ++(-1.25,0.000000) -- ++(0,1.250000) -- ++(1.25,0.000000) -- ++(0,1.250000) -- ++(-1.25,0.000000) -- ++(0,-1.250000) -- ++(-1.25,0.000000) -- ++(0,-1.250000) -- ++(1.25,0.000000) -- ++(0,-1.250000) -- ++(-1.25,0.000000) -- ++(0,1.250000) -- ++(-1.25,0.000000) -- ++(0,-1.250000) -- ++(-1.25,0.000000) -- ++(0,1.250000) -- ++(-1.25,0.000000) -- ++(0,1.250000) -- ++(1.25,0.000000) -- ++(0,1.250000) -- ++(-1.25,0.000000) -- ++(0,-1.250000) -- ++(-1.25,0.000000) -- ++(0,1.250000) -- ++(-1.25,0.000000) -- ++(0,1.250000) -- ++(1.25,0.000000) -- ++(0,-1.250000) -- ++(1.25,0.000000) -- ++(0,1.250000) -- ++(1.25,0.000000) -- ++(0,1.250000) -- ++(-1.25,0.000000) -- ++(0,1.250000) -- ++(1.25,0.000000) -- ++(0,1.250000) -- ++(-1.25,0.000000) -- ++(0,-1.250000) -- ++(-1.25,0.000000) ;

\filldraw[black] (0,0) node[left]{$A$} circle (1);
\filldraw[black] (40,0)node[above]{$B$} circle (1);
\filldraw[black] (20,0)node[left]{$C$} circle (1);
\end{scope}
\end{scope}

\begin{scope}[shift={(-100,200)}]
  \filldraw[black] (0,0) circle (1);
\filldraw[black] (-40,0) circle (1);
  \cupA{0}{0}{0}{gray}{0}{4}
  \cupA{-40}{0}{90}{magenta}{3}{1}
  \cupA{-40}{-20}{-90}{red}{3}{1}
  \cupA{-20}{-20}{-90}{magenta}{3}{1}
  \cupA{0}{0}{90}{red}{3}{1}
  \cupA{-20}{-20}{90}{red}{3}{1}
  \cupA{-20}{-20}{90}{red}{3}{1}
  \foreach \i/\j in {0/0,-20/-20}{
    \cupA{-15}{-10}{180}{black}{0.45}{0.25}
    \cupA{-20}{-10}{180}{black}{0.45}{0.25}
  \cupA{-15}{-10}{0}{blue}{0.1}{0.25}
    \begin{scope}[shift={(\i,\j)}]
  \cupA{0}{-10}{0}{blue}{0.1}{0.25}
  \cupA{-5}{-10}{0}{blue}{0.1}{0.25}
    \cupA{-5}{-10}{180}{black}{0.45}{0.25}
    \cupA{-10}{-10}{180}{black}{0.45}{0.25}
          \cupA{-10}{-10}{0}{blue}{0.1}{0.25}
    \end{scope}
  }

    \cupA{-40}{-20}{180}{black}{0.1}{0.25}
    \cupA{-35}{-25}{180}{black}{0.1}{0.25}
          \cupA{-35}{-25}{0}{blue}{0.1}{0.25}            

          \foreach \i/\j in {-32.5/-30,-37.5/-25}{
            \begin{scope}[shift={(\i,\j)}]                
\begin{scope}
          \begin{scope}[scale={0.25},rotate={-90}]
          \foreach \i in {5,10,15}{
            \foreach \j in {0,180}{
  \cupA{-\i}{0}{\j}{black}{0.1}{0.25}
            }
          }
          \cupA{0}{0}{0}{black}{0.1}{0.25}
            \cupA{-20}{0}{180}{black}{0.1}{0.25}          
          \end{scope}
\end{scope}
\end{scope}
}

 \end{scope} 
\end{tikzpicture}
\caption{$(0001)^nLLLR$ for $n=0,1,2$.  The lower boundary from $A$ to $B$
is $Y$.}
\label{0001several}
\end{figure}
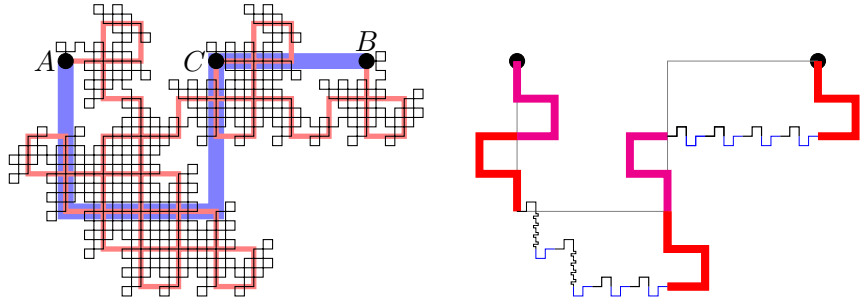

To visualise this case, see Figure~\ref{power3of0001}, which shows $(0001)^3L$.
Figure~\ref{00010001rlrrlrrrlllrrrlr} shows $(0001)^2$ applied to a configuration
of lines corresponding to $RLRRLRRRLLLRRRLR$.
The action of $(0001)$ on one edge is shown in Figure~\ref{0001A}.
This is two copies of $001$ applied to the edge, since the first application of
operation $0$ results in two edges, then $001$ is applied to both.
So we will use the same strategy as for $(001)^\infty$.
We take the shape $LLLR$ again.   Several applications of $0001$ to this
configuration are shown in Figure~\ref{0001several}, left.
Note that the upper side of this configuration
is given by four scaled copies of the lower side.  So to compute the
fractal dimension of the boundary, we just need to work this out for the
lower side, that is, whatever appears on the lower side between $A$ and $B$.
I will call this boundary $Y$.
In Figure~\ref{0001several}, right copies of $LLLR$ are shown round the
boundary of $LLLR$.  Note that care must be taken in that the image of the
point $C$ in these copies is not a boundary point, but images of
$A$ and $B$ are.
In this figure, the lines lie on three different grids, corresponding to the
blue, red, or black lines.  Imagine in each grid of lines, an edge
is coloured or left blank.  A vertex (points where grid lines intersect)
is in the boundary of $Y$ if exactly two
of the four edges are coloured.
From this figure, we see that
$$Y = \frac{1}{4}Y^{\oplus 5} + \frac{1}{16}Y^{\oplus 16} + \frac{1}{64}
Y^{\oplus 16}.$$
So the sim-value is $s\approx 1.4476$, where $\lambda = 2^{-2s}$ is a root of
\begin{equation}
  \label{eqn:001case}
  16\lambda^3 + 16\lambda^2 + 5\lambda = 1.
  \end{equation}

\section{Algorithm to find sim value}

We have found the sim value for the $(001)^\infty$ and $(0001)^\infty$
in an ad hoc way.  Now we present a more systematic method to
find the sim value for the case $B^\infty$ for any finite sequence
$B$ of $0$s and $1$s representing a sequence of our operations.

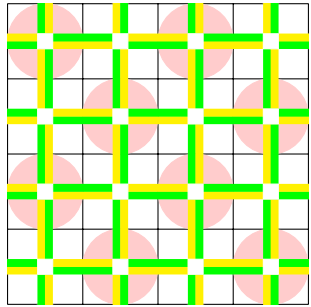
\begin{figure}
  \begin{tikzpicture}[scale = 1]
    \draw(0,0) grid (4,4);
    \begin{scope}
      \clip (0,0) rectangle (4,4);
    \foreach \i in {0,...,6}{
      \foreach \j in {0,...,5}{
        
      \filldraw[red!20!white] (\i-\j+0.5+2,\j+\i+0.5-3) circle (0.5);
      \draw[yellow, line width=3](\i-\j+0.5+1.95,\j+\i+0.5-3) -- (\i-\j+0.5+1.95,\j+\i+0.5-4);
\draw[green, line width = 3](\i-\j+0.5+2.05,\j+\i+0.5-3) -- (\i-\j+0.5+2.05,\j+\i+0.5-4);      
      \draw[green, line width=3](\i-\j+0.5+1.95,\j+\i+0.5-3) -- (\i-\j+0.5+1.95,\j+\i+0.5-2);
\draw[yellow, line width = 3](\i-\j+0.5+2.05,\j+\i+0.5-3) -- (\i-\j+0.5+2.05,\j+\i+0.5-2);      
      \draw[yellow, line width=3](\i-\j+0.5+2,\j+\i+0.5-3-0.05) -- (1+\i-\j+0.5+2,\j+\i+0.5-3-0.05);
      \draw[green, line width = 3](\i-\j+0.5+2,\j+\i+0.5-3+0.05) -- (1+\i-\j+0.5+2,\j+\i+0.5-3+0.05);
      \draw[green, line width=3](\i-\j+0.5+2,\j+\i+0.5-3-0.05) -- (-1+\i-\j+0.5+2,\j+\i+0.5-3-0.05);
\draw[yellow, line width = 3](\i-\j+0.5+2,\j+\i+0.5-3+0.05) -- (-1+\i-\j+0.5+2,\j+\i+0.5-3+0.05);      
    }}
       \end{scope}
        \foreach \i in {0,...,5}{
          \foreach \j in {0,...,5}{
\draw[white,fill=white](\i + 0.4,\j+0.4) rectangle(\i+0.6,\j+0.6);
            }}

  \end{tikzpicture}
  \caption{Colouring of line segments according to eventual boundary of
    limiting space filling curves.}
  \label{grid}
  \end{figure}

First consider a grid of the initial lines.
Remember that in a Truchet tiling path, there is a consistent
direction of travel along paths, which we denoted by arrows in
Figure~\ref{directions}.  In Figure~\ref{grid}, the background
grid is drawn in black lines.  The paths on the tiles as described in
Figure~\ref{compare} right, are drawn in yellow and green lines,
as we did in Figure~\ref{0000A}.
We will apply the same operation $B$ to the whole grid.
An sequence such as $LhRvL\dots$ tells us which lines to
follow on a grid path, and if we just take single lines, then the
result will be tessellating tiles.  The transformation of the path
depends on whether the starting tile is odd or even, and whether the
grid line crossed is horizontal or vertical.  So I put pink circles on
the even tiles to distinguish them.  In the limit of applying $B^\infty$,
the yellow side of the line will become a boundary $U$ and the green
side a boundary $V$.  The orientation of these lines is determined by whether
or not they start at an even or odd square.

\begin{figure}
  \begin{tikzpicture}[scale = 1]
    \draw[gray](1,0) grid (7,6);
    \begin{scope}
      \clip (1,0) rectangle (7,6);
    \foreach \i in {0,...,6}{
      \foreach \j in {0,...,5}{
        
      \filldraw[red!20!white] (\i-\j+0.5+2,\j+\i+0.5-3) circle (0.5);
      \draw[yellow, line width=3](\i-\j+0.5+1.95,\j+\i+0.5-3) -- (\i-\j+0.5+1.95,\j+\i+0.5-4);
\draw[green, line width = 3](\i-\j+0.5+2.05,\j+\i+0.5-3) -- (\i-\j+0.5+2.05,\j+\i+0.5-4);      
      \draw[green, line width=3](\i-\j+0.5+1.95,\j+\i+0.5-3) -- (\i-\j+0.5+1.95,\j+\i+0.5-2);
\draw[yellow, line width = 3](\i-\j+0.5+2.05,\j+\i+0.5-3) -- (\i-\j+0.5+2.05,\j+\i+0.5-2);      
      \draw[yellow, line width=3](\i-\j+0.5+2,\j+\i+0.5-3-0.05) -- (1+\i-\j+0.5+2,\j+\i+0.5-3-0.05);
      \draw[green, line width = 3](\i-\j+0.5+2,\j+\i+0.5-3+0.05) -- (1+\i-\j+0.5+2,\j+\i+0.5-3+0.05);
      \draw[green, line width=3](\i-\j+0.5+2,\j+\i+0.5-3-0.05) -- (-1+\i-\j+0.5+2,\j+\i+0.5-3-0.05);
\draw[yellow, line width = 3](\i-\j+0.5+2,\j+\i+0.5-3+0.05) -- (-1+\i-\j+0.5+2,\j+\i+0.5-3+0.05);      
    }}
       \end{scope}
        \foreach \i in {0,...,5}{
          \foreach \j in {0,...,5}{
\draw[white,fill=white](\i + 0.4,\j+0.4) rectangle(\i+0.6,\j+0.6);
            }}
        \begin{scope}[shift = {(3.5,4.5)}]
          \begin{scope}[rotate={-90}]
\draw[black, line width=1] (0,0) -- ++(0,1) 
-- ++(-1,0) -- ++(0,1) -- ++(1,0) -- ++(0,-1) -- ++(1,0) -- ++(0,1) -- ++(1,0) -- ++(0,-1) -- ++(1,0) -- ++(0,-1) -- ++(-1,0) -- ++(0,-1) -- ++(1,0) -- ++(0,1) -- ++(1,0) ;
\filldraw[red, line width=1] (0,0) circle (0.1);
            \filldraw[black, line width=1] (4,0) circle (0.1);
                    \end{scope}
        \end{scope}

        \begin{scope}[shift={(-3,5)}]
        \draw[green, line width=3](0.05,0)--(0.05,-4);
        \draw[yellow, line width=3](-0.05,0)--(-0.05,-4);
        \draw[line width = 1](0,0)--(0,-4);
        \filldraw[red, line width=1] (0,0) circle (0.1);
        \filldraw[black, line width=1] (0,-4) circle (0.1);
        \end{scope}
        \draw[->](-2,3)--(0,3);

        \node at (-3.5,3){$U$};
        \node at (-2.5,3){$V$};

        \begin{scope}[shift={(3.5,4.5)}]
          \filldraw[cyan] (1,1) circle (0.1);
          \filldraw[cyan] (2,1) circle (0.1);
          \filldraw[cyan] (2,0) circle (0.1);
          \filldraw[cyan] (2,-1) circle (0.1);
          \filldraw[cyan] (2,-2) circle (0.1);
          \filldraw[cyan] (1,-3) circle (0.1);

          \filldraw[orange] (1,-1) circle (0.1);
          \filldraw[orange] (1,-2) circle (0.1);
          \filldraw[orange] (0,-2) circle (0.1);
          \filldraw[orange] (-1,-2) circle (0.1);
          \filldraw[orange] (-1,-3) circle (0.1);          
        \end{scope}

        \begin{scope}[shift={(9.5,4.5)}]

\draw[green,line width=3](1,1.05)--(2,1.05);
\draw[yellow,line width=3](2.05,1)--(2.05,0);
\draw[yellow,line width=3](1.95,0)--(1.95,-1);
\draw[yellow,line width=3](2.05,-1)--(2.05,-2);
\draw[green,line width=3](.95,-1)--(.95,-2);
\draw[green,line width=3](-1.05,-2)--(-1.05,-3);
\draw[yellow,line width=3](0,-2.05)--(1,-2.05);
\draw[yellow,line width=3](0,-1.95)--(-1,-1.95);

          \draw[dashed,line width =2,cyan!90!black](0,0)--(1,1);
          \draw[dashed,line width =2,cyan!90!black](2,-2)--(0,-4);
\draw[dashed,line width =2,orange!90!black](-1,-3)--(0,-4);
\draw[dashed,line width =2,orange!90!black](0,0)--(1,-1);

                  \filldraw[red] (0,0) circle (0.1);
          \filldraw[black] (0,-4) circle (0.1);                  
                  \filldraw[cyan] (1,1) circle (0.1);
          \filldraw[cyan] (2,1) circle (0.1);
          \filldraw[cyan] (2,0) circle (0.1);
          \filldraw[cyan] (2,-1) circle (0.1);
          \filldraw[cyan] (2,-2) circle (0.1);
          \filldraw[cyan] (1,-3) circle (0.1);

          \filldraw[orange] (1,-1) circle (0.1);
          \filldraw[orange] (1,-2) circle (0.1);
          \filldraw[orange] (0,-2) circle (0.1);
          \filldraw[orange] (-1,-2) circle (0.1);
          \filldraw[orange] (-1,-3) circle (0.1);          
        \end{scope}
    
  \end{tikzpicture}
  \caption{$L\rightarrow 0001L$}
  \label{examplealgorithm}
  \end{figure}
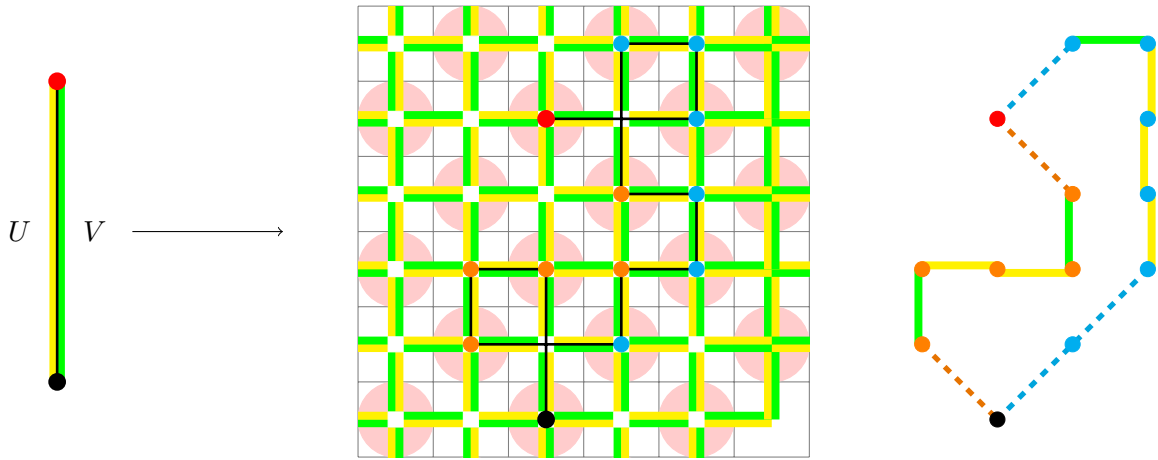

Now consider the action of our sequence of operations $B$ on a single line
segment, for example $B=0001$.  We draw the result of $0001 L$ in
Figure~\ref{examplealgorithm}.  $U$ indicates the limit of the left side
of a line $L$ going down from an even square, and $V$ is the right side,
marked in yellow and green respectively.  On this figure, the red and black dots
indicate the starting and ending points of the line respectively.  These are
both boundary points of $U$ and $V$.  Of the other vertices, those which
are in the boundary of $U$ are coloured with an orange dot, and
those which are in the boundary of $V$ are coloured with a blue dot.

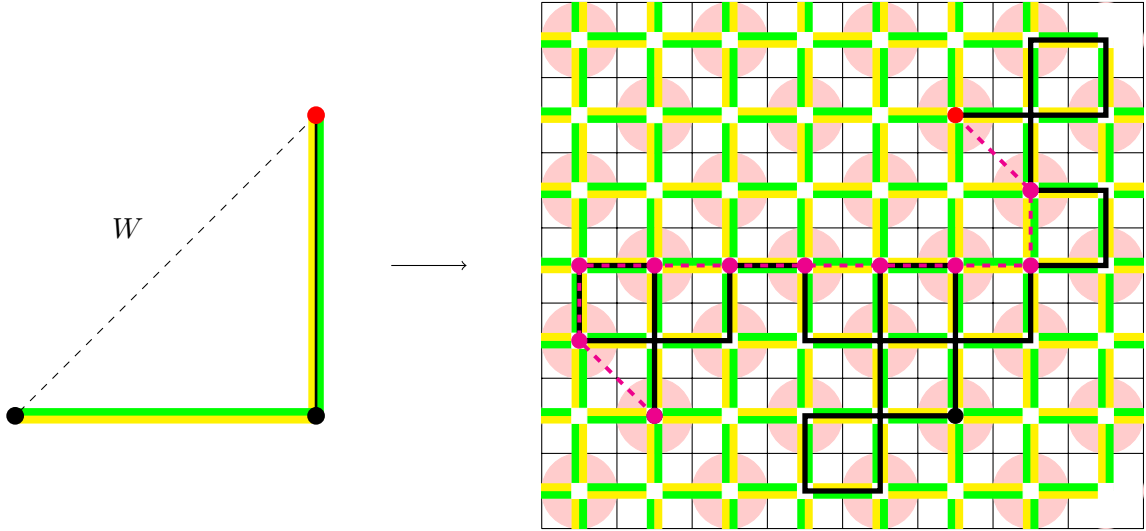
\begin{figure}
  \begin{tikzpicture}[scale = 1]
    \draw(0,1) grid (8,8);
    \begin{scope}
      \clip (0,1) rectangle (8,8);
    \foreach \i in {0,...,7}{
      \foreach \j in {0,...,7}{
        
      \filldraw[red!20!white] (2+\i-\j+0.5+2,\j+\i+0.5-3) circle (0.5);
      \draw[yellow, line width=3](\i-\j+0.5+1.95,\j+\i+0.5-3) -- (\i-\j+0.5+1.95,\j+\i+0.5-4);
\draw[green, line width = 3](\i-\j+0.5+2.05,\j+\i+0.5-3) -- (\i-\j+0.5+2.05,\j+\i+0.5-4);      
      \draw[green, line width=3](\i-\j+0.5+1.95,\j+\i+0.5-3) -- (\i-\j+0.5+1.95,\j+\i+0.5-2);
\draw[yellow, line width = 3](\i-\j+0.5+2.05,\j+\i+0.5-3) -- (\i-\j+0.5+2.05,\j+\i+0.5-2);      
      \draw[yellow, line width=3](\i-\j+0.5+2,\j+\i+0.5-3-0.05) -- (1+\i-\j+0.5+2,\j+\i+0.5-3-0.05);
      \draw[green, line width = 3](\i-\j+0.5+2,\j+\i+0.5-3+0.05) -- (1+\i-\j+0.5+2,\j+\i+0.5-3+0.05);
      \draw[green, line width=3](\i-\j+0.5+2,\j+\i+0.5-3-0.05) -- (-1+\i-\j+0.5+2,\j+\i+0.5-3-0.05);
\draw[yellow, line width = 3](\i-\j+0.5+2,\j+\i+0.5-3+0.05) -- (-1+\i-\j+0.5+2,\j+\i+0.5-3+0.05);      
    }}
       \end{scope}
        \foreach \i in {0,...,7}{
          \foreach \j in {0,...,7}{
\draw[white,fill=white](\i + 0.4,\j+0.4) rectangle(\i+0.6,\j+0.6);
        }}

        \begin{scope}[shift={(5.5,6.5)}]

\begin{scope}[rotate={-90}]
\draw[line width=2] (0,0) -- ++(0,1) 
-- ++(-1,0) -- ++(0,1) -- ++(1,0) -- ++(0,-1) -- ++(1,0) -- ++(0,1) -- ++(1,0) -- ++(0,-1) -- ++(1,0) -- ++(0,-1) -- ++(-1,0) -- ++(0,-1) -- ++(1,0) -- ++(0,1) -- ++(1,0) -- ++(0,-1) -- ++(1,0) -- ++(0,-1) -- ++(-1,0) -- ++(0,1) -- ++(-1,0) -- ++(0,-1) -- ++(-1,0) -- ++(0,-1) -- ++(1,0) -- ++(0,-1) -- ++(-1,0) -- ++(0,-1) -- ++(1,0) -- ++(0,1) -- ++(1,0) ;
\end{scope}

\draw[magenta,dashed, line width=1.5](0,0)--(1,-1)--(1,-2)
--(0,-2)--(-5,-2)--(-5,-3)--(-4,-4);

          \filldraw[red] (0,0) circle (0.1);
          \filldraw[black] (0,-4) circle (0.1);          
\filldraw[black] (-4,-4) circle (0.1);        

        \filldraw[magenta](1,-1) circle (0.1);
        \filldraw[magenta](1,-1) circle (0.1);
        \filldraw[magenta](1,-2) circle (0.1);
        \filldraw[magenta](0,-2) circle (0.1);
        \filldraw[magenta](-1,-2) circle (0.1);
        \filldraw[magenta](-2,-2) circle (0.1);
        \filldraw[magenta](-3,-2) circle (0.1);
        \filldraw[magenta](-4,-2) circle (0.1);
        \filldraw[magenta](-5,-2) circle (0.1);
        \filldraw[magenta](-5,-3) circle (0.1);        
        \filldraw[magenta](-4,-4) circle (0.1);

\end{scope}

        \begin{scope}[shift={(-3,6.5)}]          
        \draw[green, line width=3](0.05,0)--(0.05,-4);
        \draw[yellow, line width=3](-0.05,0)--(-0.05,-4);

        \draw[green, line width=3](0,-3.95)--(-4,-3.95);
        \draw[yellow, line width=3](0,-4.05)--(-4,-4.05);
        \draw[dashed](0,0)--(-4,-4);
        
        \draw[line width = 1](0,0)--(0,-4);
        \filldraw[red, line width=1] (0,0) circle (0.1);
        \filldraw[black, line width=1] (0,-4) circle (0.1);
        \filldraw[black, line width=1] (-4,-4) circle (0.1);

  \end{scope}
        \draw[->](-2,4.5)--(-1,4.5);

        \node at (-5.5,5){$W$};

  \end{tikzpicture}
  \caption{Image of $W=LR$ under $0001$}
  \label{Wcase}
  \end{figure}

\begin{figure}
\input{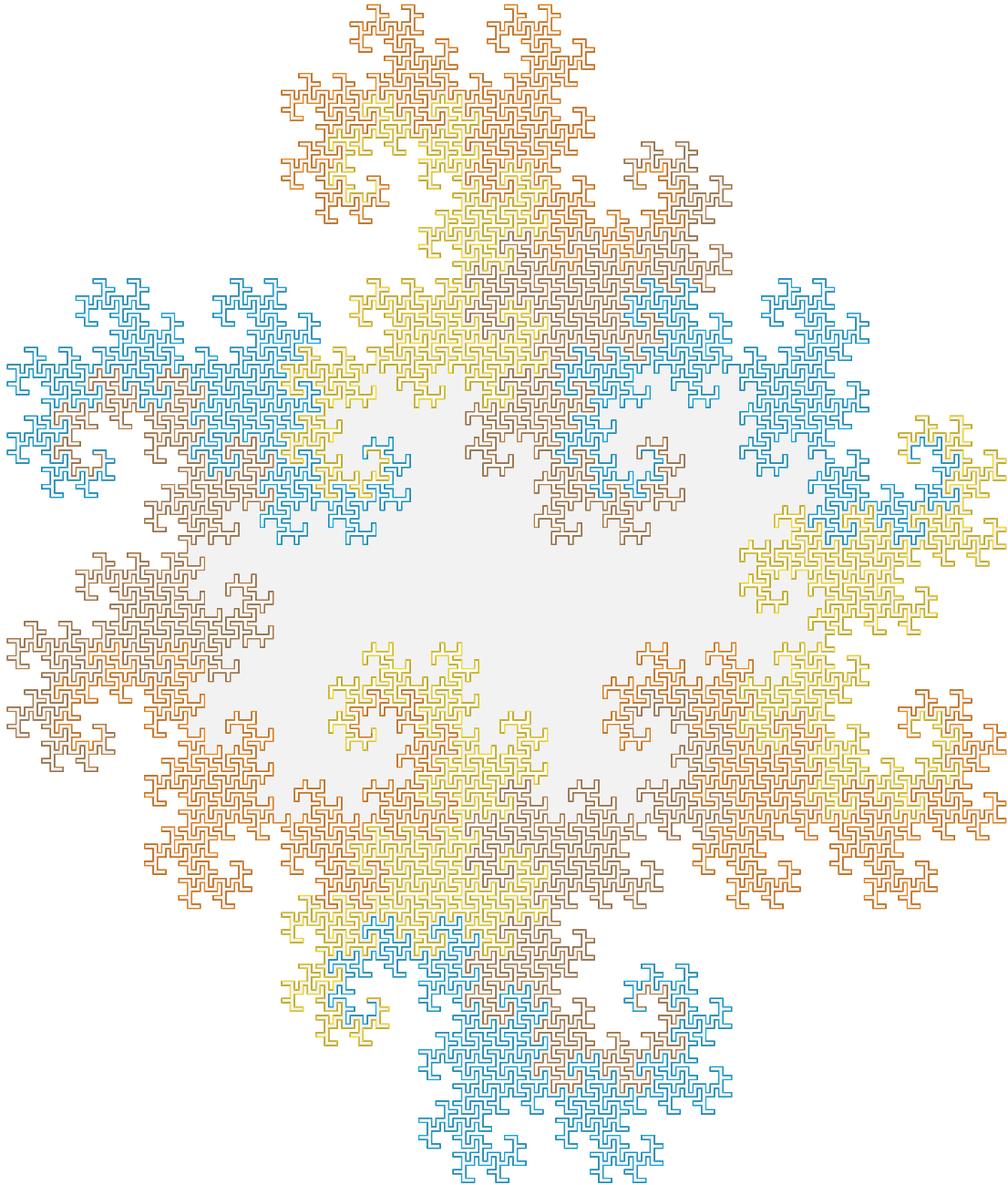}
\caption{$0000000000$ applied to $LLRLLRLLRLLR$ with diagonal line tile.}
\label{loop}
\end{figure}

\begin{figure}
  \input{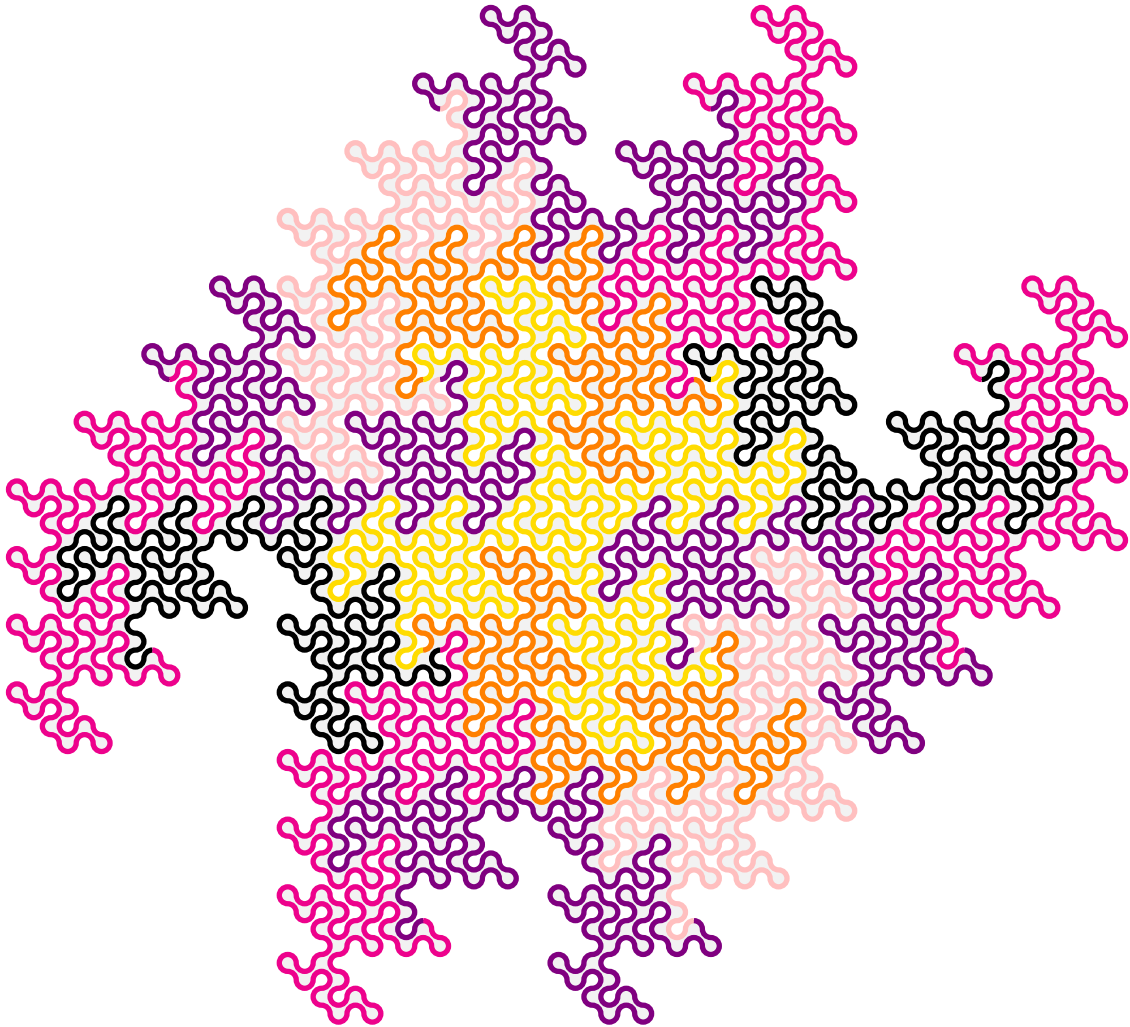}
  \caption{$00010001RLRRLRRRLLLRRRLR$ with quarter circles tile.}
  \label{00010001rlrrlrrrlllrrrlr}
  \end{figure}

\begin{figure}
  \input{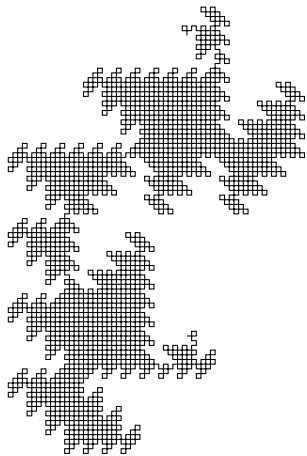}
  \caption{$(001)^3L$ with straight line tile.}
  \label{power3of0001}
  \end{figure}

Our goal is to write $U$ and $V$ in terms of smaller copies of themselves.
In Figure~\ref{examplealgorithm} we see we also need a diagonal
component, given by the limit of the concatenation of the $U$ and $V$
at right angles to each other.  We call this diagonal component $W$,
and how it transforms under $0001$ is shown in Figure~\ref{examplealgorithm}.
Let $U', V', W'$ be $U, V, W$ scaled by a factor of $\frac{1}{4}$
respectively.
Figure~\ref{examplealgorithm} right follows round the path in the
middle grid to show how $U$ and $V$ are given
in terms of $U', V', W'$.  The situation for $W$ is shown in
Figure~\ref{Wcase}.  The magenta points are in the boundary of $W$, and following the
dashed line path through them we can read off $W=W'+U'+6V'+W'$.
From these figures we have the system of equations:
\begin{alignat*}{4}
U & {}={} &  2U' & {}+{} & 2V' & {}+{} & 2W' \\
V & {}={} &  3U' & {}+{} & V' & {}+{} & 3W' \\
W & {}={} &  6U' & {}+{} & 2V' & {}+{} & 2W'
\end{alignat*}

We can rewrite this as a matrix equation
\begin{equation}
  {\mathbf U} = M  {\mathbf U}'
  \;\text{ where }\;
  M = \left(
  \begin{array}{ccc}
    2&2&2\\
    3&1&3\\
    6&2&2
    \end{array}
  \right),
  \;
  {\mathbf U} = (U,V,W)^{T}.
  \end{equation}
Here ${\mathbf U}' = (U',V',W')=(sU,sV,sW)$ where the scale factor
is $s=\frac{1}{4}$.  In general, if $B$ has length $n$,
the factor $s$ will be $2^{-n/2}$,
To find the fractal dimension of $U$ and $V$, we want to
find the sim value of some shape $X$ composed of a union of
$U, V, W$.
Suppose we have $X=\alpha U + \beta V + \gamma W$, which
we write as
$$X = A\mathbf U,$$
where $A = (\alpha,\beta,\gamma)$.
In order to write $X$ in terms of $X', X''$ and $X'''$
(which are equal to $sX, s^2X, s^3X$), we write
everything in terms of $X'''$.  By repeated application of $M$,
we see that
\begin{alignat*}{2}  
  X & {}={} &  AM^3\mathbf U'''\\
  X' & {}={} &  AM^2\mathbf U'''\\
  X'' & {}={} &  AM\mathbf U'''\\    
  X''' & {}={} &  A\mathbf U'''
  \end{alignat*}
Writing $X = (a,b,c)(X',X'',X''')$ for some non negative integers
$a,b,c$ implies that
$$A(M^3- aM^2 - bM - cI)=0.$$
For this to be possible, we must find non negative
integers $a,b,c$ with
$\det(M^3- aM^2 - bM - cI)=0.$
To do this, it is sufficient to find
the characteristic polynomial $c_M$ of the matrix $M$.  Provided
all but the leading term have non positive coefficients,
our $a, b, c$ are provided, and we can just take $X=U$.
In this example, it turns out that $c_M(x)=x^3 - 5x^2 - 16x -16$.
These give the coefficients in
Equation~\ref{eqn:001case}, which determine the sim value as
previously computed.

To summarize, the outline of the algorithm to compute
the sim value of the boundary of $B^\infty L$ (where $L$
is an initial line)  for
some finite binary string $B$ of length one is as follows:
\begin{enumerate}
\item Draw $BU, BV$ and $BW$, as in examples in Figures~\ref{examplealgorithm} and
  \ref{Wcase}.
\item Find a matrix $M$ to express $U, V, W$ in terms
  of $sU, sV, sW$ for $s=2^{-n/2}$.
\item If possible,
  find non negative integers not all zero $\alpha, \beta, \gamma$ and
  non negative not all zero $a,b,c$ such
  that if $c_M(x)=x^3 - ax^2 - bx - c$ then
   $(\alpha, \beta, \gamma)c_M(M)=0$.
If this is possible, we can proceed to the next step.
\item Let $\lambda$ be a real root of $c_M$.  Then the sim value of
  $B^\infty L$ is
  $$s=\frac{-2\log(\lambda)}{n\log(2)}.$$
  \end{enumerate}

 \section{Future work}
It remains to show that the algorithm for finding the sim value
gives the fractal dimension, for which we need an algorithm for finding
a Moran open set.  We also need to automate the process of computing the
matrix $M$, by writing a program for computing the boundaries of
$BL$ and $BLR$, which was done by hand in Figures~\ref{examplealgorithm}
and \ref{Wcase}.  Then it would be possible to investigate the
possible values and the
upper bound on the fractal dimensions obtained from this process.  Many of the
diagrams for this paper were produced by a JavaScript program found on the
explanation page of \cite{Verrill}, which can be referred to to reproduce the latex figures.
This program is a starting point for further computations.
Future work would be to investigate fractals produced from other hinged
tilings, or with other related replacement rules.

    \bibliographystyle{alpha}
\bibliography{ms}

\end{document}